\documentclass[12pt,reqno]{amsart}

\usepackage{geometry,color}
\usepackage{hyperref}

\DeclareMathOperator{\csch}{csch}

\usepackage{graphicx}
\usepackage{amssymb,amsmath}
\usepackage{amsthm}
\usepackage[pdf]{pstricks}
\usepackage{empheq}
\usepackage{pgfplots}

\usepackage{epsfig}
\usepackage{extarrows}
\usepackage{amsfonts}
\usepackage{caption}
\usepackage{tikz}

\newtheorem{theorem}{Theorem}[section]

\newtheorem{proposition}[theorem]{Proposition}
\newtheorem{corollary}[theorem]{Corollary}

\numberwithin{equation}{section}

\theoremstyle{definition}

\title[Stokes waves in a finite-depth fluid]
{Stokes waves in a finite-depth fluid \\ via a conformal transformation}

\date{\today}

\author[A. Aghayev]{Ali Aghayev}
\address[A. Aghayev]{Department of Mathematics and Statistics, McMaster University, Hamilton, Ontario, Canada, L8S 4K1}
\email{aghayeva@mcmaster.ca}

\author[E. Byrnes]{Eleanor Byrnes}
\address[E. Byrnes]{Department of Applied Mathematics, University of Washington, Seattle, WA USA 98195-3925}
\email{elbyrnes@uw.edu}

\author[D. E. Pelinovsky]{Dmitry E. Pelinovsky}
\address[D. E. Pelinovsky]{Department of Mathematics and Statistics, McMaster University, Hamilton, Ontario, Canada, L8S 4K1}
\email{pelinod@mcmaster.ca}

\begin{document}
\maketitle

\begin{abstract}
Euler equations for a spatially periodic surface in a two-dimensional fluid of a finite depth are reduced to a system of evolution equations after a conformal transformation. Compared to the previously developed formal 
approach, we show that the system of evolution equations is complemented 
by two constraints which determine the new depth in 
conformal variables and the Bernoulli integration constant. A locally unique continuation of the traveling Stokes waves with respect to the wave speed is studied by using arguments based on the implicit function theorem. 
We show that the constraints on the traveling wave profile recover the correct small-amplitude expansion of Stokes waves in physical variables. We also prove for the co-periodic stability problem that the real unstable eigenvalues may bifurcate from the zero eigenvalue only at a critical point of the wave energy. Finally, we perform numerical approximations of the wave profiles for both small and large amplitudes and illustrate different bifurcations along the solution family. 
\end{abstract}


\section{Introduction}

Surface gravity waves in a fluid of finite depth are essential for understanding laboratory experiments in water tanks and dynamics of ocean flows. Unlike in a deep fluid, the finite depth modifies both the dispersion and nonlinear interactions of water waves and complicates the analysis of the underlying equations. We are particularly interested in the stability of the traveling periodic waves called Stokes waves since understanding the linear stability of traveling periodic waves is important for predicting wave breaking and the formation of rogue waves \cite{KPS08}.

In deep water, the linear stability of Stokes waves has been studied by Longuet-Higgins \cite{L-H 1,L-H 2}, who suggested that the change of stability occurs when the wave speed attains an extremum with respect to wave steepness. Numerical computations by Tanaka \cite{Tanaka83,Tanaka85} showed that this criterion is not accurate and that the exchange of stability occurs at the extremal point of the wave energy. Saffman \cite{Saffman85} provided an  analytical proof of these numerical findings based on Zakharov’s equations  by showing that the zero eigenvalue of the co-periodic stability problem has a higher algebraic multiplicity at the critical point of the wave energy. In the recent work \cite{DP25}, this property was rigorously proven by using the conformal transformation of the fluid domain to a vertical semi-infinite strip. It was also shown in \cite{DMP26} that the Stokes wave at the critical point of the wave speed is associated with the zero eigenvalue of a higher geometric multiplicity and the same algebraic multiplicity. The conformal transformation of Stokes waves in a deep fluid within the so-called Dyachenko--Zakharov (DZ) formulation \cite{DKSZ96,Zakharov96} was also used to study the existence of Stokes waves \cite{BDT2000a,BDT2000b}, 
the well-posedness of Euler's equations \cite{Ifrim2}, 
and the modulational stability of Stokes waves \cite{DDS24,DMP26,DS23,DS23SAPM}.

Far fewer analytical and numerical studies exist for Stokes waves in a finite-depth fluid. Stability of solitary waves was investigated numerically 
in \cite{Tanaka86} and the stability criterion based on the extremal points of the wave energy was extended in \cite{Z-S86} to both solitary and periodic waves based on Zakharov's equations. The conformal transformation in the DZ formulation was used for the analysis of well-posedness of Euler's equations \cite{Ifrim1}, the analysis of the existence of traveling solitary waves \cite{Mats} (see also earlier works \cite{Buffoni1,Buffoni2} in the case of a nonzero surface tension), the numerical study of steep wave profiles \cite{BDS26}, and the derivation of a new evolution equation combining smooth and peaked profiles \cite{LP25,NPW25}.

The purpose of this work is to address the system of evolution equations in the periodic fluid domain of finite depth after a conformal transformation. Compared to the previous works \cite{DKSZ96,Zakharov96}, we show that the system of evolution equations in the DZ formulation is complemented by two constraints which were previously ignored. The first of these constraints was previously discovered in \cite{DH19} and used to compute Stokes waves in \cite{BDS26}, but the second, connecting the surface elevation and the mean level of the fluid to the Bernoulli constant has not yet been reported in the literature. The two constraints determine the new depth in conformal variables and the Bernoulli integration constant. The derivation of the system of constrained evolution equations is reported in Section \ref{sec-equations}, where we also elaborate on the mass, momentum, and energy conservation in Propositions \ref{prop-conserved} and \ref{prop-energy}.

We study the full system of constrained evolution equations and obtain 
the following four main results of this study.
\begin{enumerate}
\item In Section \ref{sec-existence}, we formulate the conditions for the locally unique continuation of the wave profile with respect to the wave speed, see Proposition \ref{prop-2-parameters}. The horizontal momentum and the energy of Stokes waves are characterized in Proposition \ref{prop-momentum} and Corollary \ref{prop-momentum-1-par}.

\medskip

\item In Section \ref{sec-stability}, we prove that the zero eigenvalue 
in the co-periodic stability problem admits a higher algebraic multiplicity at each critical point of the horizontal momentum, see Theorem \ref{theorem-multiplicity} and Corollary \ref{cor-multiplicity}. We also show that each critical point of the horizontal momentum coincides with a critical point of the wave energy. The main tool of the theory is  the mass, momentum, and energy conservation in the linearized system of equations, as described in Propositions \ref{prop-lin-energy} and \ref{prop-orth}, as well as  
computations of the geometric and algebraic multiplicities of the zero eigenvalue. 

\medskip

\item In Section \ref{sec-small}, we show that the small-amplitude expansions for the constrained evolution equations 
in conformal variables correctly recover the Stokes expansion of the wave profile in physical variables. These computations agree with similar studies in other formulations of Euler's equations such as the Craig--Sulem (CS) formulation used in \cite{BMV23} and the Ablowitz--Fokas--Musslimani (AFM) formulation used in \cite{CD23,CDT22}. We also establish the correspondence between parameters of traveling waves used in the DZ formulation \cite{DKSZ96,Zakharov96} and parameters of traveling waves used in the AFM formulation \cite{CD23,CDT22}, see Proposition \ref{prop-comparison}.

\medskip

\item In Section \ref{sec:numerics}, we extend the numerical results of \cite{BDS26,SB25} to the traveling waves under the two constraints. We give plots of the Stokes wave profiles, as well as their derivatives with respect to the wave speed and the conformal depth. Furthermore, we demonstrate that the wave speed, energy, and Bernoulli constant satisfy the asymptotic theory for steep Stokes waves from \cite{LHF1,LHF2}. Within this study, we provide numerical evidence  that the critical points of the speed no longer correspond to an increase in the dimension of the kernel of the system's Jacobian, as in the case of infinite depth \cite{DMP26}. Instead, this corresponds to the critical points of the bifurcation parameter $\gamma$ for a fixed conformal depth, which precede critical points of the bifurcation parameter $\gamma$ for a fixed physical depth. Critical points of energy occur before the two critical points of $\gamma$, whereas critical points of the speed occur after the two critical points of $\gamma$. We conjecture that this sequence of critical points is repeated infinitely many times when the family of steep Stokes waves is extended to the limiting wave with the peaked profile.
\end{enumerate}

Section \ref{sec-conclusion} concludes the paper with the summary and further directions.

\medskip

{\bf Acknowledgements.} The authors thank B. Deconinck and S. Dyachenko for many stimulating discussions. 
A. A. acknowledges the partial support of the MSc fellowship at McMaster University. 
E. B. acknowledges the partial support of the Wan fellowship.
D. P. acknowledges the partial support of the NSERC Discovery grant.

\section{Derivation of a system of constrained evolution equations}
\label{sec-equations}

Here, we transform Euler's equations for the spatially periodic surface of a two-dimensional fluid of a finite depth to the system of constrained evolution equations by using a conformal transformation. In what follows, the space of square--integrable and $2\pi$-periodic functions is denoted by $L^2_{\rm per}$ and equipped with the inner product $\langle \cdot, \cdot \rangle$ and the induced norm $\| \cdot \|_{L^2_{\rm per}}$. The Sobolev space of functions in $L^2_{\rm per}$ with the first derivative being also square--integrable is denoted by $H^1_{\rm per}$ and equipped with the norm $\| \cdot \|_{H^1_{\rm per}}$.

\subsection{Conformal mapping of the fluid domain}

We consider a two--dimensional, spatially $2\pi$-periodic domain for a fluid of a finite depth $h_0 > 0$. The free surface is located at $y=\eta(x,t)$. The velocity potential $\varphi = \varphi(x,y,t)$ satisfies the Laplace equation in the fluid domain
\[
\mathcal{D}_{h_0}(t) = \{ (x,y) : \quad x \in \mathbb{T}, \quad  -h_0 \le y \le \eta(x,t)\},
\]
under the Neumann boundary condition $\varphi_y |_{y = -h_0} = 0$ 
and the periodic boundary conditions on the circle $\mathbb{T}$ of length $2\pi$. Kinematic and dynamic boundary conditions hold at the free surface, according to Euler's equations
\begin{align}
y = \eta(x,t) : \quad 
\begin{cases}
\eta_t + \varphi_x \eta_x - \varphi_y = 0,  \\
\varphi_t + \frac{1}{2} (\varphi_x)^2 + \frac{1}{2} (\varphi_y)^2 + \eta = 0, 
\end{cases}
\label{euler-1-2}
\end{align}
where the gravitational constant $g$ is normalized to unity and the atmospheric pressure $p_0$ is set to zero.

Let $z = z(u,t)$ be a conformal transformation of the rectangular region 
\[
\tilde{\mathcal{D}}_h(t) = \{ u \in \mathbb{C} : \quad {\rm Re}(u) \in \mathbb{T}, \quad  -h(t) \le {\rm Im}(u) \le 0\}
\]
into the fluid domain $\mathcal{D}_{h_0}(t)$, where $h = h(t)$ is the new depth in conformal variables and  ${\rm Im}(u) = 0$ corresponds to the parameterization of the free surface as $\xi + i \eta = z(u,t)$, where $(\xi,\eta)$ are now functions of $(u,t)$ with $u \in \mathbb{T}$. The other dynamical variable is the potential $\psi = \varphi(x,y,t)$ evaluated at the free surface $y = \eta(x,t)$. After the conformal transformation, $\psi$ is a function of $(u,t)$ with $u \in \mathbb{T}$. 

As shown in \cite[Appendix A]{LP25}, the mapping between $\xi$ and $\eta$ at ${\rm Im}(u) = 0$ is given by
\begin{equation}
\label{eq:xuKh}
\eta = h - h_0 + T_h (\xi - u),
\end{equation}
where $\eta$ and $\xi - u$ are periodic functions on $\mathbb{T}$ and where the nonlocal operator $T_h$ is 
defined by its Fourier multiplicative symbol
\[
\widehat{(T_h)}_n = i\tanh(nh), \qquad n\in\mathbb Z.
\]
The linear operator $T_h : L^2_{\rm per} \to L^2_{\rm per}$ is bounded and skew-adjoint, with the kernel ${\rm Ker}(T_h) = {\rm span}(1)$. The inverse operator $T_h^{-1}$ is only defined on the periodic functions in $L^2_{\rm per}$ with zero mean, denoted by $L^2_{\rm per} |_{\{1\}^{\perp}}$. The uniquely defined inverse operator $T_h^{-1} : L^2_{\rm per} |_{\{1\}^{\perp}} \to L^2_{\rm per} |_{\{1\}^{\perp}}$ has its Fourier multiplicative symbol
\[
\widehat{(T_h^{-1})}_n = \left\{ \begin{array}{cl} -i\coth(nh), &\qquad n \in\mathbb Z \backslash \{0\}, \\
0, & \qquad n = 0. \end{array} \right.
\]

Integrating (\ref{eq:xuKh}) over the period on $\mathbb{T}$ yields the following local relation between the physical depth $h_0 > 0$ and the conformal depth $h = h(t)$, 
\begin{equation}\label{eq:DepthConstraint}
h = h_0 + \bar{\eta}, \quad \mbox{\rm where} \;\; \bar{\eta} = \bar{\eta}(t) := \frac{1}{2\pi} \int_{\mathbb{T}} \eta(u,t)\,du.
\end{equation}
The conformal depth $h$ depends on time $t$ if and only if the mean value 
of the surface elevation $\bar{\eta}$ in the conformal variable $u \in \mathbb{T}$ depends on time $t$. Since $\bar{\eta}$ is not generally conserved in time, the time dependence $h = h(t)$ is generally nontrivial.

Inverting (\ref{eq:xuKh}) with $T_h^{-1} : L^2_{\rm per} |_{\{1\}^{\perp}} \to L^2_{\rm per} |_{\{1\}^{\perp}}$ yields 
\begin{equation}
\label{xi-expressions}
\xi - u = T_h^{-1} \eta,
\end{equation}
where the constant of integration for the mean value of the periodic function 
$\xi - u$ can be set to zero due to the invariance of the periodic domain $\mathbb{T}$ with respect to the choice of the left boundary point. It follows from 
(\ref{xi-expressions}) and the chain rule that 
\begin{equation}
\label{xi-expressions-der}
\xi_u = 1 + K_h \eta, \quad \xi_t = T_h^{-1} \eta_t + h'(t) (\partial_h T_h^{-1}) \eta,
\end{equation}
where $K_h = T_h^{-1}\partial_u : {\rm Dom}(K_h) \subset L^2_{\rm per} \to L^2_{\rm per}$ is a self-adjoint, positive semi-definite, linear operator with the domain ${\rm Dom}(K_h) = H^1_{\rm per}$, defined by its Fourier multiplicative symbol
\begin{equation}
\label{eq:KhDefinition}
(K_h)_n = \left\{ \begin{array}{ll} n \coth(nh), \quad & n\in\mathbb{Z} \backslash \{0\}, \\
0, \quad  & n = 0. \end{array} \right.
\end{equation}
The linear operator $K_h$ admits the kernel ${\rm Ker}(K_h) = {\rm span}(1)$.

\subsection{Derivation of equations of motion}

As shown in \cite[Appendix A]{LP25}, see also \cite[Appendix A]{DLK16} for the deep fluid, equations of motion 
in conformal variables are derived from the Lagrangian functional
\begin{align*}
\mathcal{L}(\psi,\xi,\eta,h,f)
&=
\int_{\mathbb{T}}
\psi \bigl(\eta_t \xi_u - \eta_u \xi_t\bigr)\,du
-
\frac12
\int_{\mathbb{T}}
\left( \eta^2 \xi_u - \psi \,T_h\psi_u \right)\,du \\
&\quad
+
\int_{\mathbb{T}}
f\bigl(\eta - h + h_0 - T_h(\xi - u)\bigr)\,du.
\end{align*}
The first term of the Lagrangian represents the kinetic energy of the fluid, 
the second term represents the negative potential energy of the fluid, 
and the last term contains a Lagrange multipliers $f$  introduced to 
preserve the constraint (\ref{eq:xuKh}). By writing $g = 2\pi \int_{\mathbb{T}} f du$, 
we rewrite the Lagrangian in the equivalent form
\begin{align*}
\mathcal{L}(\psi,\xi,\eta,h,f,g)
&=
\int_{\mathbb{T}}
\psi \bigl(\eta_t \xi_u - \eta_u \xi_t\bigr)\,du
-
\frac12
\int_{\mathbb{T}}
\left( \eta^2 \xi_u - \psi \,T_h\psi_u \right)\,du \\
&\quad
+
\int_{\mathbb{T}}
f\bigl(\eta - h + h_0 - T_h(\xi - u)\bigr)\,du
+
g \left( \frac{1}{2\pi} \int_{\mathbb{T}} \eta du  - h +h_0 \right),
\end{align*}
where the Lagrange multiplier $g$ preserves the constraint (\ref{eq:DepthConstraint}) explicitly and 
we now impose the zero-mean constraint on $f$.

Variations of the Lagrangian $\mathcal{L}$ in $f$ and $g$ recover the constraints (\ref{eq:xuKh}) and (\ref{eq:DepthConstraint}). Variations of the Lagrangian $\mathcal{L}$ in $\psi$, $\eta$, $\xi$, and $h$ yield the following system of equations:
\begin{align}
\delta_{\psi}: & \quad 
\eta_t \xi_u-\eta_u \xi_t+T_h \psi_u=0, \label{eq-1} \\
\delta_{\eta}: & \quad 
\psi_u \xi_t-\psi_t \xi_u - \eta \xi_u + f + \frac{1}{2\pi}g = 0, \label{eq-2}\\
\delta_{\xi}: & \quad 
\psi_t\, \eta_u - \psi_u\, \eta_t +\eta\,\eta_u + T_h f = 0, \label{eq-3}\\
\delta_{h}: & \quad 
\frac12 \langle \psi, (\partial_h T_h)\,\psi_u \rangle
- \langle f, (\partial_h T_h)(\xi-u) \rangle - g = 0,\label{eq-4}
\end{align}
where we have used the zero-mean constraint on $f$. Note that the last equation (\ref{eq-4}) was not discussed in the previous works \cite{DKSZ96,Zakharov96}.

We can now express Lagrange multipliers $f$ and $g$ from (\ref{eq-3}) and (\ref{eq-4}):
\begin{equation}
\label{f-expression}
f = -T_h^{-1} (\psi_t\, \eta_u - \psi_u\, \eta_t +\eta\,\eta_u)
\end{equation}
and 
\begin{equation}
\label{g-expression}
g = \frac12 \langle \psi, (\partial_h T_h)\,\psi_u \rangle 
- \langle f,(\partial_h T_h)(\xi-u) \rangle.
\end{equation}
Substituting these expressions into (\ref{eq-2}) closes the system of equations 
(\ref{eq-1}) and (\ref{eq-2}) for $(\psi,\eta,\xi,h)$ augmented by the constraints (\ref{eq:xuKh}) and (\ref{eq:DepthConstraint}). In particular, 
the variable $\xi$ and its derivatives in $u$ and $t$ can be eliminated from (\ref{xi-expressions}) and (\ref{xi-expressions-der}). 

\subsection{Conserved quantities}

To proceed further, we introduce Bernoulli's integration constant 
$B = B(t)$ for the $2\pi$-periodic function $\psi(\cdot,t)$ by writing
\begin{equation}
\label{Bernoulli}
\psi(u,t) = \int_0^t B(t) dt + \tilde{\psi}(u,t),
\end{equation}
where $\tilde{\psi}(\cdot,t)$ is required to have the zero mean value for any $t \in \mathbb{R}$. The conserved quantities of the mass, 
vertical momentum, and horizontal momentum of the fluid were obtained in physical coordinates \cite{BO82}. The following proposition 
recovers the same quantities from the constrained system of equations. 

\begin{proposition}
	\label{prop-conserved}
The system of equations (\ref{eq-1}), (\ref{eq-2}), and (\ref{eq-3}) augmented by the decomposition (\ref{Bernoulli}) admits conservation of the following three integral quantities:
\begin{align}
M(\eta) &=
\langle \eta,\xi_u \rangle = \langle \eta, (1 + K_h \eta) \rangle,
\label{mass-M} \\
Q(\tilde{\psi},\eta) &= 
\langle \tilde{\psi}, \xi_u \rangle = \langle \tilde{\psi}, (1+K_h\eta) \rangle, 
\label{momentum-Q} \\
P(\tilde{\psi},\eta) &=
- \langle \tilde{\psi}, \eta_u \rangle,
\label{momentum-P}
\end{align} 
if and only if $B(t)$, $M(\eta)$, and $g$ are related by 
\begin{equation}
\label{g-relation}
g = 2 \pi B(t) + M(\eta).
\end{equation}
\end{proposition}

\begin{proof}
By taking mean values of 
(\ref{eq-1}), (\ref{eq-2}), and (\ref{eq-3}), we obtain 
\begin{align*}
\frac{d}{dt} M(\eta) = 0, \quad 
\frac{d}{dt} Q(\tilde{\psi},\eta) = 0, \quad 
\frac{d}{dt} P(\tilde{\psi},\eta) = 0, 
\end{align*} 
where conservation of $Q$ in (\ref{momentum-Q}) is justified if and only if 
the constraint (\ref{g-relation}) is satisfied. 
\end{proof}

Combining (\ref{xi-expressions}), (\ref{f-expression}), (\ref{g-expression}), (\ref{Bernoulli}), and (\ref{g-relation}), we obtain the integral constraint in the form:
\begin{equation}
\label{constraint-preliminary}
2 \pi B(t) + M(\eta) = \frac12 \langle \tilde{\psi}, (\partial_h T_h)\,\tilde{\psi}_u \rangle
+ \langle (\tilde{\psi}_t\, \eta_u - \tilde{\psi}_u\, \eta_t + B(t) \eta_u +  \eta \eta_u), (\partial_h T_h^{-1}) \eta \rangle, 
\end{equation}
where we have used $\partial_h T_h^{-1} = -T_h^{-1} (\partial_h T_h) T_h^{-1}$ and the skew-adjointness of $T_h^{-1}$.

The mass integral $M(\eta)$ is conserved in time by Proposition \ref{prop-conserved}. In what follows, we set the mass integral to $0$ for the proper normalization of the surface elevation, $M(\eta) = 0$.
The constraint $M(\eta) = 0$ implies 
the zero mean value of the surface elevation $\eta$ in the physical coordinate $x$. This condition defines the fluid depth $h_0$ uniquely. 

Substituting $f$ and $g$ from (\ref{f-expression}) and (\ref{g-relation}) with $M(\eta) = 0$ into (\ref{eq-2}) and substituting $\psi$ from (\ref{Bernoulli}) with $B = B(t)$ into (\ref{eq-1}) and (\ref{eq-2}) yields a system of two coupled evolution equations:
\begin{align}
\eta_t \xi_u -\eta_u \xi_t+ T_h \tilde{\psi}_u =& 0, 
\label{time-1} \\
\tilde{\psi}_u \xi_t - \tilde{\psi}_t \xi_u - (B + \eta) \xi_u -T_h^{-1} (\tilde{\psi}_t\, \eta_u - \tilde{\psi}_u\, \eta_t + (B + \eta) \,\eta_u) + B =& 0, \label{time-2}
\end{align}
where $\xi_u$ and $\xi_t$ can be expressed from (\ref{xi-expressions-der}). 
The system of coupled equations (\ref{time-1}) and  (\ref{time-2}) is constrained by the local constraint 
(\ref{eq:DepthConstraint}) and the nonlocal constraint (\ref{constraint-preliminary}), which is rewritten with $M(\eta) = 0$ as 
\begin{equation}
\label{time-constraint}
2 \pi B = \frac{1}{2} \langle \tilde{\psi}, (\partial_h T_h)\,\tilde{\psi}_u \rangle 
+ \langle (\tilde{\psi}_t\, \eta_u - \tilde{\psi}_u\, \eta_t + (B+\eta)\,\eta_u), (\partial_h T_h^{-1}) \eta \rangle.
\end{equation}
The following proposition establishes the energy conservation.

\begin{proposition}
	\label{prop-energy}
	The energy $H(\tilde{\psi},\eta)$ given by 
	\begin{equation}
	\label{energy-H}
	H(\tilde{\psi},\eta) = \langle (\eta^2 + 2 B \eta), (1 + K_h \eta) \rangle - \langle \tilde{\psi}, T_h \partial_u \tilde{\psi} \rangle
	\end{equation}
	is conserved in time if and only if the 
	integral constraint (\ref{time-constraint}) holds.
\end{proposition}

\begin{proof}
	To derive the expression for $\frac{d}{dt} H(\tilde{\psi},\eta)$, we multiply (\ref{time-1}) by $\tilde{\psi}_t$, 
	multiply (\ref{time-2}) by $\eta_t$, integrate both with respect to $u$ over the period, 
	and add them together. We get 
	\begin{align}
	& \int_{\mathbb{T}} \left[ \tilde{\psi}_t T_h \partial_u \tilde{\psi} + (\eta_t \tilde{\psi}_u - \tilde{\psi}_t \eta_u) \xi_t 
	- (B+\eta) \eta_t \xi_u - B \eta_t K_h \eta - \frac{1}{2} \eta_t K_h \eta^2 \right] du \notag \\
	& \quad + \int_{\mathbb{T}} \eta_t T_h^{-1} (\eta_t \tilde{\psi}_u - \tilde{\psi}_t \eta_u) du + 2 \pi B h'(t) = 0,
    \label{tech-energy-1}
	\end{align}
where we have used 
	\begin{align*}
	h'(t) = \frac{1}{2\pi} \int_{\mathbb{T}} \eta_t du,
	\end{align*}
which follows from (\ref{eq:DepthConstraint}). Differentiating (\ref{energy-H}) in $t$ and using $M(\eta) = 0$, we write 
\begin{align*}
\frac{1}{2} \frac{d}{dt} H(\tilde{\psi},\eta) &= 
\langle \eta_t, \eta (1 + K_h \eta) \rangle + \frac{1}{2} \langle \eta_t, K_h \eta^2 \rangle + B \langle \eta_t, (1 + 2 K_h \eta) \rangle 
- \langle \tilde{\psi}_t, T_h \partial_u \tilde{\psi} \rangle \\
& \quad + \frac{1}{2} \langle (\eta^2 + 2 B \eta), (\partial_h K_h) \eta \rangle - \frac{1}{2} \langle \tilde{\psi}, (\partial_h T_h) \partial_u \tilde{\psi} \rangle,
\end{align*}
where we have used self-adjointness of $T_h \partial_u$ and $K_h$. Substituting (\ref{xi-expressions-der}) into  (\ref{tech-energy-1}) yields
	\begin{align*}
	\frac{1}{2} \frac{d}{dt} H(\tilde{\psi},\eta) &= h'(t) \left[ 2\pi B + \langle (\eta_t \tilde{\psi}_u - \tilde{\psi}_t \eta_u), (\partial_h T_h^{-1}) \eta \rangle \right]  \\
	& \quad + \frac{1}{2} h'(t) \left[ \langle (\eta^2 + 2 B \eta), (\partial_h K_h) \eta \rangle - \langle \tilde{\psi}, (\partial_h T_h)  \tilde{\psi}_u \rangle \right],
	\end{align*}
where we have used skew-adjointness of $T_h^{-1}$. Since $K_h = T_h^{-1} \partial_u$, integration by parts yields that $\frac{d}{dt} H(\tilde{\psi},\eta) = 0$ if and only if the integral constraint (\ref{time-constraint}) is satisfied for any given $h'(t)$.
\end{proof}

\section{Existence of traveling waves}
\label{sec-existence}

Here, we derive Babenko's equation, named after \cite{Babenko87}, for the profile of the traveling wave 
solutions of the constrained evolution equations (\ref{time-1}), (\ref{time-2}), and (\ref{time-constraint}). Then, we analyze a locally 
unique continuation of the traveling wave profile with respect to the parameter 
of the wave speed $c$ after fixing the fluid depth $h_0 > 0$. Parameters 
of the conformal depth $h$ and the Bernoulli constant $B$ are uniquely defined 
in terms of the wave speed $c$ for the traveling wave solutions.

\subsection{Babenko's equation for the finite depth}

Steady traveling waves with constant speed $c$ are expressed by 
the $2\pi$-periodic profiles 
\begin{equation}
\label{TW-profiles}
\eta(u,t) = \eta(u-ct), \qquad
\tilde{\psi}(u,t) = \tilde{\psi}(u-ct).
\end{equation}
The conformal depth $h$ in the transformed domain $\tilde{\mathcal{D}}_h$ and the Bernoulli constant $B$ are independent of time $t$ for the steady traveling waves. It follows from (\ref{xi-expressions-der}) for the traveling waves that 
\begin{equation}
\label{xi-trav}
\xi_u = 1 + K_h \eta, \quad \xi_t = -c K_h \eta.
\end{equation}
Substituting (\ref{TW-profiles}) into (\ref{time-1}) and integrating with respect to $u$ yields 
\begin{equation}
\label{psi-trav}
\tilde{\psi} = c\,T_h^{-1}\eta + \tilde{\psi}_0,
\end{equation}
where the constant term $\tilde{\psi}_0$ can be set to $0$ without loss of generality. Substituting (\ref{xi-trav}) and (\ref{psi-trav}) into (\ref{time-2}) yields  Babenko's equation:
\begin{equation}
\label{eq:Babenko}
(c^{2} - 2B) K_h \eta - \eta
= \frac12 K_h \eta^{2} + \eta\,K_h\eta.
\end{equation}
The only variable of (\ref{eq:Babenko}) is the traveling wave coordinate $u-ct$, which we denote again by $u$. The nonlocal constraint (\ref{time-constraint}) can be rewritten in the equivalent form:
\begin{equation}\label{eq:DepthIdentity}
B = \frac{1}{4\pi} \langle  (c^{2} - 2B) \eta - \eta^2, (\partial_h K_h) \eta \rangle,
\end{equation}
where we have again used that $(\partial_h T_h^{-1}) = -T_h^{-1} (\partial_h T_h) T_h^{-1}$, as well as the skew-adjointness of $T_h^{-1}$ and $\partial_u$.
Although the system of equations (\ref{eq:Babenko}) and (\ref{eq:DepthIdentity}) 
depends on three parameters $c$, $h$, and $B$, we recall the existence of  the additional local constraint
\begin{equation}
\label{eq:Local}
h = h_0 + \bar{\eta}, \quad \mbox{\rm where} \quad \bar{\eta} = \frac{1}{2\pi} \int_{\mathbb{T}} \eta(u) du,
\end{equation}
which follows from (\ref{eq:DepthConstraint}).

The system of equations \eqref{eq:Babenko}, (\ref{eq:DepthIdentity}), and (\ref{eq:Local}) admits a variational formulation with the action functional
\begin{align}
\Lambda(\eta,h) &= \frac{1}{2} \int_{\mathbb{T}}
\left[ \eta\,(c^{2}K_h - 1)\eta - \eta^{2}K_h\eta - 2 B \eta (1 + K_h \eta) \right]\,du \notag \\
& \quad + g \left( \frac{1}{2\pi} \int_{\mathbb{T}} \eta du - h + h_0 \right),
\label{eq:ActionFunctional}
\end{align}
where $g$ is the Lagrange multiplier introduced to impose the local constraint (\ref{eq:Local}). Since $c$ and $B$ are fixed parameters, 
variation of $\Lambda(\eta,h)$ in $(\eta,h)$ yields the system of equations:
\begin{align*}
\delta_{\eta}: & \quad 
(c^2 K_h - 1) \eta - \frac{1}{2} K_h \eta^2 - \eta K_h \eta - B - 2 B K_h \eta + \frac{g}{2\pi} = 0, \\
\delta_{h}: & \quad 
\frac{1}{2} \langle (c^2 - 2B) \eta - \eta^2, (\partial_h K_h) \eta \rangle - g = 0.
\end{align*}
Taking the mean value of the first equation yields 
$g = 2\pi B$. Eliminating $g = 2\pi B$ from the two Euler--Lagrange equations yields the system \eqref{eq:Babenko} and (\ref{eq:DepthIdentity}).

\subsection{Characterization of the traveling waves}

We consider smooth solutions of Babenko's equation (\ref{eq:Babenko}) 
with the profile $\eta \in C^{\infty}_{\rm per}$ and the Bernoulli constant $B$ under the constraints (\ref{eq:DepthIdentity}) and (\ref{eq:Local}). 
To resolve the constraints, we introduce the bifurcation parameter $\gamma > 0$ by $\gamma^2 = c^2 - 2B > 0$ and rewrite (\ref{eq:Babenko}) in the explicit form:
\begin{equation}
\label{BE}
\gamma^2 K_h \eta - \eta = \frac{1}{2} K_h \eta^2 + \eta K_h \eta.
\end{equation}
Similarly to the case of $h = \infty$ in \cite{BDT2000a,BDT2000b}, 
the local and global bifurcation theories imply that there exists a smooth solution $\eta \in C^{\infty}_{\rm per}$ of the normalized Babenko's equation (\ref{BE}) in an open region of the $(\gamma,h)$-parameter plane, that is, $\eta = \eta(u;\gamma,h)$. Substituting the profile $\eta \in C^{\infty}_{\rm per}$ into the integral constraint (\ref{eq:DepthIdentity}) yields the unique expression 
for $B = B(\gamma,h)$ given by 
\begin{equation}
\label{IC}
B = \frac{1}{4\pi} \langle \gamma^2 \eta - \eta^2, (\partial_h K_h)\eta \rangle.
\end{equation}
The value of $h = h(\gamma)$ is to be defined from the only remaining constraint (\ref{eq:Local}) imposed for a fixed value of the physical depth $h_0$. We assume that there exists a solution of (\ref{eq:Local}) for $(\gamma,h) = (\gamma_*,h_*)$. We also assume 
that the solution satisfies the non-degeneracy assumption 
${\rm Ker}(\mathcal{L}_h) = {\rm span}(\eta')$, where 
$\eta'$ stands for the derivative of $\eta \in C^{\infty}_{\rm per}$ with respect to $u$ and $\mathcal{L}_h : H^1_{\rm per} \subset L^2_{\rm per} \to L^2_{\rm per}$ is the linearized operator for the nonlinear equation (\ref{BE}) given by 
$$
\mathcal{L}_h = \gamma^2 K_h - 1 - K_h \eta - \eta K_h - K_h(\eta \; \cdot).
$$
The following result guarantees a locally unique continuation of the solution $(\eta,B)$ with respect to the wave speed $c$ for a fixed value of the physical depth $h_0$.

\begin{proposition}
	\label{prop-2-parameters}
	Let $(\eta,B) \in C^{\infty}_{\rm per} \times \mathbb{R}$ be a solution of the system (\ref{BE}) and (\ref{IC}) subject to the constraint $h = h_0 + \bar{\eta}$ with fixed $h_0 > 0$ for one point $(\gamma,h) = (\gamma_*,h_*)$. 
	Assume that ${\rm Ker}(\mathcal{L}_h) = {\rm span}(\eta')$ 
	and define $c_* := \sqrt{\gamma_*^2 + 2B}$. If the following constraint is satisfied for $(\gamma,h) = (\gamma_*,h_*)$, 
	\begin{equation}
	\label{smooth-cont}
	1 - \overline{\partial_h \eta} \neq 0, \quad \overline{\partial_h \eta} = \frac{1}{2\pi} \int_{\mathbb{T}} \partial_h \eta du,
	\end{equation}
then the solution $(\eta,B) \in C^{\infty}_{\rm per} \times \mathbb{R}$ of the system (\ref{eq:Babenko}) and (\ref{eq:DepthIdentity}) is uniquely and smoothly continued with respect to $(c,h)$ near $(c_*,h_*)$. Moreover, there exists a unique smooth solution of $h = h_0 + \bar{\eta}$ with fixed $h_0 > 0$ for $(c,h)$ near $(c_*,h_*)$. In addition, we have $\partial_c B \neq c$.
\end{proposition}

\begin{proof}
	Since the system (\ref{BE}) and (\ref{IC}) is smooth in $(\gamma,h)$, the implicit function theorem gives a unique smooth continuation of the solution pair $(\eta,B) \in C^{\infty}_{\rm per} \times \mathbb{R}$ in $(\gamma,h)$ near $(\gamma_*,h_*)$ under the assumption ${\rm Ker}(\mathcal{L}_h) = {\rm span}(\eta')$. Furthermore, since the constraint $h = h_0 + \bar{\eta}$ is satisfied at $(\gamma,h) = (\gamma_*,h_*)$ and 
	the condition (\ref{smooth-cont}) is fulfilled at $(\gamma,h) = (\gamma_*,h_*)$, the implicit function theorem gives a unique smooth solution of 
    the constraint
\begin{equation}
\label{constr-again}
h_0 = h - \frac{1}{2\pi} \int_{\mathbb{T}} \eta(u;\gamma,h) du,
\end{equation}	
for $(\gamma,h)$ near $(\gamma_*,h_*)$. It remains to transfer these results to a continuation in $c$ near $c_* := \sqrt{\gamma_*^2 + 2B(\gamma_*,h_*)}$. 
	
By the implicit function theorem for $c^2 = \gamma^2 + 2B(\gamma,h)$, 
there exists a unique smooth continuation of the solution $(\eta,B) \in C^{\infty}_{\rm per} \times \mathbb{R}$  in $(c,h)$ near $(c_*,h_*)$ if 
the following constraint is satisfied at $(\gamma,h) = (\gamma_*,h_*)$,
\begin{equation}
\label{B-der-IFT}
\gamma + \partial_{\gamma} B(\gamma,h) \neq 0.
\end{equation}
We show that (\ref{B-der-IFT}) is equivalent to (\ref{smooth-cont}). Taking partial derivatives of (\ref{BE}) in $\gamma$ and $h$ under the assumption ${\rm Ker}(\mathcal{L}_h) = {\rm span}(\eta')$, we get 
\begin{align*}
\begin{cases}
    \mathcal{L}_h \partial_{\gamma} \eta = -2 \gamma K_h \eta, \\
    \mathcal{L}_h \partial_{h} \eta = -\gamma^2 (\partial_h K_h) \eta + \frac{1}{2}  (\partial_h K_h) \eta^2 + \eta  (\partial_h K_h) \eta.
    \end{cases}
\end{align*}
Since all functions belong to $H^1_{\rm per}$ and $\mathcal{L}_h$ is self-adjoint in $L^2_{\rm per}$, we obtain 
\begin{align*}
\begin{cases}
    \langle \mathcal{L}_h \partial_{\gamma} \eta, \partial_h \eta \rangle = -2 \gamma \langle K_h \eta, \partial_h \eta \rangle,\\
    \langle \mathcal{L}_h \partial_{h} \eta, \partial_{\gamma} \eta \rangle = -\gamma^2 \langle (\partial_h K_h) \eta, \partial_{\gamma} \eta \rangle  
    + \frac{1}{2}  \langle (\partial_h K_h) \eta^2,\partial_{\gamma} \eta \rangle + \langle \eta  (\partial_h K_h) \eta, \partial_{\gamma} \eta \rangle,
        \end{cases}
\end{align*}
with $\langle \mathcal{L}_h \partial_{\gamma} \eta, \partial_h \eta \rangle = \langle \mathcal{L}_h \partial_{h} \eta, \partial_{\gamma} \eta \rangle$. 
Differentiating (\ref{IC}) in $\gamma$ yields
\begin{align*}
2\pi \partial_{\gamma} B &= \gamma \langle \eta, (\partial_h K_h) \eta \rangle + \frac{1}{2} \langle (\partial_h K_h) \eta, (\gamma^2 - 2\eta) \partial_{\gamma} \eta \rangle  + \frac{1}{2}  \langle (\gamma^2 \eta - \eta^2), (\partial_h K_h) \partial_{\gamma} \eta \rangle \\
    &= \gamma \langle \eta, (\partial_h K_h) \eta \rangle - \langle \mathcal{L}_h \partial_{h} \eta, \partial_{\gamma} \eta \rangle \\
&= \gamma \langle \eta, (\partial_h K_h) \eta \rangle + 2 \gamma \langle K_h  \eta, \partial_h \eta \rangle \\
&= \gamma \partial_h \langle \eta, K_h \eta \rangle,
\end{align*}
where we used the fact that $\partial_h K$ is self-adjoint in $L^2_{\rm per}$. Recalling the mass constraint 
$$
M(\eta) = \langle 1, \eta \rangle + \langle \eta, K_h \eta \rangle = 0,
$$
we finally obtain
$$
\gamma + \partial_{\gamma} B(\gamma,h) = \gamma \left[ 1 - \frac{1}{2\pi} \partial_h \langle 1, \eta \rangle \right],
$$
and since $\gamma \neq 0$, the constraint (\ref{B-der-IFT}) is equivalent to (\ref{smooth-cont}). Thus, the solution $(\eta,B) \in C^{\infty}_{\rm per} \times \mathbb{R}$ of the system (\ref{eq:Babenko}) and (\ref{eq:DepthIdentity})  is uniquely and smoothly continued in $(c,h)$ near $(c_*,h_*)$. Furthermore, the solution $(\gamma,h)$ of the constraint (\ref{constr-again}) near $(\gamma_*,h_*)$ transfers to the solution 
$(c,h)$ of the same constraint (\ref{constr-again}) near $(c_*,h_*)$.
Finally, we use the chain rule for $c^2 = \gamma^2 + 2B(\gamma,h)$ and obtain 
$$
c - \partial_c B = c - \frac{c \partial_{\gamma} B}{\gamma + \partial_{\gamma} B} = \frac{c \gamma}{\gamma + \partial_{\gamma} B} \neq 0,
$$
since $c, \gamma \neq 0$. Hence, we get $\partial_c B \neq c$. 
\end{proof}

\subsection{Derivative equations}

Let $(\eta,B) \in C^{\infty}_{\rm per} \times \mathbb{R}$ be a solution of the system 	(\ref{eq:Babenko}) and (\ref{eq:DepthIdentity}) defined for $(c,h)$ near $(c_*,h_*)$ by Proposition \ref{prop-2-parameters} under the assumptions ${\rm Ker}(\mathcal{L}_h) = {\rm span}(\eta')$ 
and $1 - \overline{\partial_h \eta} \neq 0$. To proceed further, we overload notations and denote the solutions by $\eta = \eta(u;c,h)$ and $B = B(c,h)$. Let us also overload notation for 
$\mathcal{L}_h : H^1_{\rm per} \subset L^2_{\rm per} \to L^2_{\rm per}$ as the linearized operator for the nonlinear equation (\ref{eq:Babenko}):
\begin{equation}
    \label{def-L}
\mathcal{L}_h = (c^2 -2 B) K_h - 1  - K_h \eta - \eta K_h - K_h (\eta \; \cdot ).
\end{equation}
We also denote Babenko's equation (\ref{eq:Babenko}) as $S_h(\eta) = 0$ with 
\begin{equation}
    \label{def-S}
S_h(\eta) = (c^2 - 2B) K_h \eta - \eta - \eta K_h \eta - \frac{1}{2} K_h \eta^2,
\end{equation}
and define 
\begin{equation}
    \label{def-S-der}
\partial_h S_h(\eta) = (c^2 - 2B) (\partial_h K_h) \eta - \eta (\partial_h K_h) \eta - \frac{1}{2} (\partial_h K_h) \eta^2.
\end{equation}
Assuming ${\rm Ker}(\mathcal{L}_h) = {\rm span}(\eta')$ 
and $1 - \overline{\partial_h \eta} \neq 0$, we are allowed by Proposition \ref{prop-2-parameters} to differentiate (\ref{eq:Babenko}) with $B = B(c,h)$ in $u$, $h$, and $c$. This yields three linear equations:
\begin{align}
\mathcal{L}_h \eta' &= 0, \label{derivative-1} \\
\mathcal{L}_h \partial_h \eta &= -\partial_h S_h(\eta) + 2 (\partial_h B) K_h \eta, \label{derivative-2} \\
\mathcal{L}_h \partial_c \eta &= -2 c K_h \eta + 2 (\partial_c B) K_h \eta, \label{derivative-3}
\end{align}
where $\eta'$ stands for the derivative of $\eta = \eta(u;c,h)$ with respect to $u$. Since $\partial_c B \neq c$ by Proposition \ref{prop-2-parameters}, the inhomogeneous term in (\ref{derivative-3}) is not identically equal to $0$. 

Similarly, assuming ${\rm Ker}(\mathcal{L}_h) = {\rm span}(\eta')$ 
and $1 - \overline{\partial_h \eta} \neq 0$, we are allowed by Proposition \ref{prop-2-parameters} to differentiate the integral constraint (\ref{eq:DepthIdentity}) in $h$ and $c$. This yields two integral equations:
\begin{equation}
\label{derivative-Ber-1}
2\pi  \left[ 1 + \frac{1}{2\pi} \langle \eta, (\partial_h K_h) \eta \rangle \right] \partial_h B = \langle \partial_h \eta, \partial_h S_h(\eta) \rangle 
+ \frac{1}{2}  \langle (c^2 - 2B) \eta - \eta^2, (\partial_h^2 K_h)\, \eta \rangle 
\end{equation}
and 
\begin{equation}
\label{derivative-Ber-2}
2\pi  \left[ 1 + \frac{1}{2\pi} \langle \eta, (\partial_h K_h) \eta \rangle \right] \partial_c B = \langle \partial_c \eta, \partial_h S_h(\eta) \rangle 
+ c \langle \eta, (\partial_h K_h)\, \eta \rangle. 
\end{equation}
The derivative equations are used in the stability analysis of the traveling waves. 

\subsection{The horizontal momentum and energy of the traveling wave}

We compute the horizontal momentum (\ref{momentum-P}) and the energy 
(\ref{energy-H}) at the traveling wave with profile $\eta \in C^{\infty}_{\rm per}$, for which $\tilde{\psi} = c T_h^{-1} \eta$ 
and obtain
\begin{equation}
\label{hor-momentum}
\mathcal{P}(c,h) = c \langle K_h \eta, \eta \rangle
\end{equation}
and
\begin{equation}
\label{energy-wave}
\mathcal{H}(c,h) = c^2 \langle K_h \eta, \eta \rangle + \langle \eta^2, (1+K_h \eta) \rangle,
\end{equation}
since $M(\eta) = 0$. The following proposition computes the partial derivative of $\mathcal{P}(c,h)$ in $(c,h)$.

\begin{proposition}
	\label{prop-momentum}
Under the assumptions of Proposition \ref{prop-2-parameters}, we have
	\begin{equation}
	\label{derivative-momentum}
	\frac{\partial \mathcal{P}}{\partial c} = \langle K_h \eta, \eta \rangle - c \langle 1, \partial_c \eta \rangle, \quad 
	\frac{\partial \mathcal{P}}{\partial h} = -c \langle 1, \partial_h \eta \rangle.
	\end{equation}
\end{proposition}

\begin{proof}
	By using  the smooth continuation 
	of solutions in $(c,h)$ by Proposition \ref{prop-2-parameters}, 
	we obtain by differentiating $M(\eta) = 0$ in $c$ and $h$ that 
	\begin{align*}
	\frac{\partial}{\partial c} M(\eta) = \langle (1+2K_h \eta), \partial_c \eta \rangle = 0
	\end{align*}
	and
	\begin{align*}
	\frac{\partial}{\partial h} M(\eta) = \langle (1+2K_h \eta), \partial_h \eta \rangle +  \langle \eta, (\partial_h K_h) \eta \rangle = 0.
	\end{align*}
	Then, we obtain 
	\begin{align*}
	\frac{\partial \mathcal{P}}{\partial c} &=	\langle K_h \eta, \eta \rangle + 2 c \langle K_h \eta, \partial_c \eta \rangle \\
	&= \langle K_h \eta, \eta \rangle + c \langle (1 + 2K_h \eta), \partial_c \eta \rangle - c \langle 1, \partial_c \eta \rangle \\
	&= \langle K_h \eta, \eta \rangle - c \langle 1, \partial_c \eta \rangle
	\end{align*}
	and 
	\begin{align*}
	\frac{\partial \mathcal{P}}{\partial h} &=	2 c \langle K_h \eta, \partial_h \eta \rangle + c \langle \eta, (\partial_h K_h) \eta \rangle \\
	&= c \langle (1 + 2K_h \eta), \partial_h \eta \rangle + c \langle \eta, (\partial_h K_h) \eta \rangle  - c \langle 1, \partial_h \eta \rangle \\
	&= - c \langle 1, \partial_h \eta \rangle,
	\end{align*}
	which yields the result.
\end{proof}

\begin{corollary}
	\label{prop-momentum-1-par}
Under the assumptions of Proposition \ref{prop-2-parameters}, let us define $h = h(c)$ as the solution of 
the constraint $h = h_0 + \bar{\eta}$ with fixed $h_0 > 0$. Then, we have
	\begin{equation}
\label{derivative-momentum-1-par}
\frac{d}{dc} \mathcal{P}(c,h(c)) = \langle K_h \eta, \eta \rangle - 2\pi c h'(c)
\end{equation}	
and 
\begin{equation}
\label{der-energy} 
\frac{d}{dc} \mathcal{H}(c,h(c)) = 2c \frac{d}{dc} \mathcal{P}(c,h(c)).
\end{equation}	
\end{corollary}

\begin{proof}
By differentiating (\ref{eq:Local}) in $c$, we obtain 
	\begin{equation}
	\label{derivative-constraints}
	(1 - \overline{\partial_h \eta}) h'(c) - \overline{\partial_c \eta} = 0,
	\end{equation}
To show that $1 - \overline{\partial_h \eta} \neq 0$, we emphasize that $\partial_h \eta$  in (\ref{derivative-constraints}) is obtained from $\eta = \eta(u;c,h)$, whereas $\partial_h \eta$ in the constraint (\ref{smooth-cont}) is obtained from $\eta = \eta(u;\gamma,h)$ due to overloaded notations. However, due to the invertibility of the mapping $c^2 = \gamma^2 + 2B$ in Proposition \ref{prop-2-parameters}, the unique smooth continuation of solutions of $h = h_0 + \bar{\eta}$ as $h = h(\gamma)$ in $\gamma$ implies the unique smooth continuation of solutions of $h = h_0 + \bar{\eta}$ as 
$h = h(c)$ in $c$, so that $1 - \overline{\partial_h \eta} \neq 0$ in (\ref{derivative-constraints}) follows from the constraint (\ref{smooth-cont}). 

By using (\ref{derivative-momentum}) and (\ref{derivative-constraints}), 
we get  
\begin{align*}
\frac{d}{dc} \mathcal{P}(c,h(c))  &= \frac{\partial \mathcal{P}}{\partial c} + h'(c) \frac{\partial \mathcal{P}}{\partial h} \\
&= \langle K_h \eta, \eta \rangle - c \langle 1, \partial_c \eta + h'(c) \partial_h \eta \rangle \\
&= \langle K_h \eta, \eta \rangle - 2\pi c h'(c),
\end{align*}
which yields (\ref{derivative-momentum-1-par}). To show (\ref{der-energy}), 
we use the action functional $\Lambda$ given by (\ref{eq:ActionFunctional}) with $g = 2\pi B$ 
and $M(\eta) = 0$ to obtain 
\begin{align*}
\frac{\partial \Lambda}{\partial c} = c \langle K_h \eta, \eta \rangle 
+ \langle \partial_c \eta, \underbrace{(c^2 K_h - 1) \eta - \eta K_h \eta - \frac{1}{2} K_h \eta^2 - 2B K_h \eta}_{=0} \rangle
\end{align*}
and
\begin{align*}
\frac{\partial \Lambda}{\partial h} = \langle \partial_c \eta, \underbrace{(c^2 K_h - 1) \eta - \eta K_h \eta - \frac{1}{2} K_h \eta^2 - 2B K_h \eta}_{=0} \rangle + \frac{1}{2} \underbrace{\langle (c^2 \eta - \eta^2 - 2B \eta), (\partial_h K_h) \eta \rangle}_{=4\pi B} - 2 \pi B, 
\end{align*}
where we have used both Babenko's equation (\ref{eq:Babenko}) and 
the constraint (\ref{eq:DepthIdentity}). Thus, we obtain 
$\frac{d \Lambda}{dc} = \frac{\partial \Lambda}{\partial c} = \mathcal{P}(c,h(c))$. Since 
$\Lambda = c \mathcal{P}(c,h(c)) - \frac{1}{2} \mathcal{H}(c,h(c))$, 
differentiating $\Lambda$ in $c$ yields (\ref{der-energy}).
\end{proof}

\section{Stability of traveling waves}
\label{sec-stability}

Here, we reformulate the time-dependent system in the traveling frame, perform linearization around a traveling wave profile, and derive the spectral stability problem for co-periodic perturbations. Then, we analyze the geometric and algebraic multiplicity of the zero eigenvalue to prove that a higher algebraic multiplicity of the zero eigenvalue occurs at each critical point of the horizontal momentum. By Corollary \ref{prop-momentum-1-par}, particularly, (\ref{der-energy}), every critical point of the horizontal momentum coincides with a critical point of the wave energy.

\subsection{Reformulation of the system in the traveling frame}

Passing to the traveling wave coordinate with the transformation
\begin{equation}
\label{time-TW}
\eta(u,t) = \eta(u-ct,t), \qquad
\tilde{\psi}(u,t) = \tilde{\psi}(u-ct,t),
\end{equation}
and assuming that $h = h(t)$, we obtain from (\ref{xi-expressions-der}) that 
\begin{equation}
\label{xi-derivative}
\xi_u = 1 + K_h \eta, \quad \xi_t = -c K_h \eta + T_h^{-1} \eta_t + h'(t) (\partial_h T_h^{-1}) \eta.
\end{equation}
Furthermore, written in the traveling frame as in (\ref{time-TW}), we transform $\tilde{\psi}$ by introducing a new function $\zeta$, 
\begin{equation}
\label{time-psi}
\tilde{\psi} = c\,T_h^{-1}\eta + \zeta, 
\end{equation}
which yields 
\begin{equation}
\label{psi-derivative}
\tilde{\psi}_u = c K_h \eta + \zeta_u, \quad \tilde{\psi}_t = c\,T_h^{-1}\eta_t + \zeta_t + c h'(t) (\partial_h T_h^{-1}) \eta.
\end{equation}
Equations (\ref{time-1}) and (\ref{time-2}) after the transformations (\ref{time-TW}), (\ref{xi-derivative}), and (\ref{psi-derivative}) are rewritten in the explicit form:
\begin{align}
(1 + K_h \eta) \eta_t - \eta_u T_h^{-1} \eta_t  = -T_h \zeta_u + h'(t) \eta_u (\partial_h T_h^{-1}) \eta, & \label{time-3} \\
(1 + K_h \eta) \zeta_t + T_h^{-1} (\eta_u \zeta_t) + 2c T_h^{-1} \eta_t - \zeta_u (T_h^{-1} \eta_t) - T_h^{-1} (\zeta_u \eta_t) &=  S_h(\eta) \notag\\
 - c h'(t) (\partial_h T_h^{-1}) \eta + h'(t) \zeta_u (\partial_h T_h^{-1}) \eta, & \label{time-4}
\end{align}
where $S_h(\eta) = 0$ expresses Babenko's equation (\ref{eq:Babenko}) as in (\ref{def-S}). After expansion, integration by parts, and using (\ref{time-3}), the constraint (\ref{time-constraint}) is rewritten in the form:
\begin{align}
2 \pi B &= \frac{1}{2} \langle  (c^2 - 2B) \eta - \eta^2, (\partial_h K_h)\,\eta \rangle \notag\\
& \quad + \frac{1}{2} \langle \zeta, (\partial_h T_h)\,\zeta_u \rangle
+ \langle (\zeta_t \eta_u - \zeta_u \eta_t + c\eta_t), (\partial_h T_h^{-1})\,\eta \rangle.
\label{time-constraint-expanded}
\end{align}
The mass conservation $M(\eta)$ is still written in the form (\ref{mass-M}), but the other three conserved quantities (\ref{momentum-Q}), (\ref{momentum-P}), and (\ref{energy-H}) are rewritten in variables $(\zeta,\eta)$ with the help of (\ref{time-psi}). Hence, we overload notations for $Q$, $P$, $H$ and rewrite these conserved quantities as 
\begin{align}
Q(\zeta,\eta) &= \langle \zeta, (1+K_h \eta) \rangle, 
\label{momentum-Q-again} \\
P(\zeta,\eta) &= c \langle K_h \eta, \eta \rangle - \langle  \zeta, \eta_u \rangle, \label{momentum-P-again} \\
H(\zeta,\eta) &= c^2 \langle K_h \eta, \eta \rangle + \langle \eta^2, (1 + K_h \eta) \rangle - 2c \langle \zeta,\eta_u \rangle  - \langle \zeta, T_h \partial_u \zeta \rangle,\label{energy-H-again}
\end{align}
where we have used the orthogonality of $T_h^{-1} \eta$ to $(1+K_h \eta)$, the skew-adjoint property of $T_h$, $T_h^{-1}$, $\partial_u$, and the mass constraint $M(\eta) = 0$. By Proposition \ref{prop-energy}, the integral constraint (\ref{time-constraint-expanded}) is {\em equivalent} to the conservation of $H(\zeta,\eta)$.

\subsection{Linearized equation of motion}

Let $(\eta,B) \in C^{\infty}_{\rm per} \times \mathbb{R}$ be the solution pair 
satisfying Babenko's equation (\ref{eq:Babenko}) and 
the associated constraint (\ref{eq:DepthIdentity}) by Proposition \ref{prop-2-parameters}. We substitute $\eta \to \eta + v$ and $\zeta \to w$ 
with the co-periodic perturbation $(v,w)$ into the system 
(\ref{time-3}), (\ref{time-4}), and  (\ref{time-constraint-expanded}).
We further use the relation (\ref{eq:DepthConstraint}) to write 
$h(t) = h_0 + \bar{\eta} + \bar{v}(t)$, where 
\begin{equation}
\label{mean-term}
\bar{\eta} = \frac{1}{2\pi} \int_{\mathbb{T}} \eta(u) du, \quad 
\bar{v}(t) = \frac{1}{2\pi} \int_{\mathbb{T}} v(u,t) du.
\end{equation}
Keeping only the linear terms with respect to the perturbation $(v,w)$ in the system (\ref{time-3}) and (\ref{time-4}) yields the linearized equations:
\begin{align}
\mathcal{M}_h v_t  &= -T_h \partial_u w + \bar{v}' \eta' (\partial_h T_h^{-1}) \eta,  \label{lin-1} \\
\mathcal{M}^*_h w_t + 2c T_h^{-1} v_t  &=  \mathcal{L}_h v + \bar{v} \partial_h S_h(\eta) - 2 (\delta B) K_h \eta - c \bar{v}' (\partial_h T_h^{-1}) \eta, & \label{lin-2}
\end{align}
where  $\delta B$ is the perturbation of the Bernoulli constant $B$, 
the linear operators $\mathcal{M}_h,\mathcal{M}_h^* : L^2_{\rm per} \to L^2_{\rm per}$ are given by  
\begin{align*}
\mathcal{M}_h &= (1 + K_h \eta) - \eta' T_h^{-1}, \\
\mathcal{M}_h^* &= (1 + K_h \eta) + T_h^{-1} (\eta' \cdot ), 
\end{align*}
the linearized Babenko's operator $\mathcal{L}_h : H^1_{\rm per} \subset L^2_{\rm per} \to L^2_{\rm per}$ is  defined by (\ref{def-L}), and the derivative of coefficients of Babenko's equation $\partial_h S_h(\eta)$ is defined by (\ref{def-S-der}).
Keeping only the linear terms with respect to the perturbation $(v,w)$ in the
constraint (\ref{time-constraint-expanded}) yields the linearized constraint:
\begin{align}
2 \pi (\delta B) \left[ 1 + \frac{1}{2\pi} \langle \eta, (\partial_h K_h) \eta \rangle \right] &= \langle v, \partial_h S_h(\eta) \rangle 
 + \frac{1}{2} \bar{v} \langle (c^2 - 2B) \eta - \eta^2, (\partial_h^2 K_h)\, \eta \rangle  \notag \\
& \quad + \langle  w_t, \eta' (\partial_h T_h^{-1})\,\eta \rangle + 
c \langle v_t, (\partial_h T^{-1}_h) \eta \rangle. \label{lin-constraint}
\end{align}

The following proposition relates the linearized equations of motion 
to the corrections to the mass, vertical and horizontal momenta, and energy.

\begin{proposition}
	\label{prop-lin-energy}
Corrections to the mass (\ref{mass-M}), vertical momentum (\ref{momentum-Q-again}), and horizontal momentum (\ref{momentum-P-again}) are given by 	
	\begin{align}
	\delta M &= \langle (1 + 2 K_h \eta), v \rangle + \bar{v} \langle \eta, (\partial_h K_h) \eta \rangle = 0, 	\label{lin-mass} \\
	\delta Q &= \langle (1+K_h \eta), w \rangle, 	\label{lin-momentum-vert}  \\
	\delta P &= 2c \langle K_h \eta, v \rangle - \langle \eta', w \rangle
	+ c \bar{v} \langle \eta, (\partial_h K_h) \eta \rangle, \label{lin-momentum}
	\end{align}
whereas the correction to the energy (\ref{energy-H-again}) is given by $\delta H = 2 c \delta P$. If $(v,w)$ is a suitable solution of the 
	linearized system (\ref{lin-1})--(\ref{lin-2}), then $\delta M$, $\delta P$, and $\delta Q$ are independent of time $t$.
\end{proposition}

\begin{proof}
	Validity of the expressions in (\ref{lin-mass}), (\ref{lin-momentum-vert}), and (\ref{lin-momentum}) follows by substitution $\eta \to \eta + v$ and $\zeta \to w$ into $M$, $Q$, and $P$ and the linearization. We now verify their conservation directly from the linearized equations (\ref{lin-1})--(\ref{lin-2}).
	
Projection of (\ref{lin-1}) to $1$ yields
$$
\frac{d}{dt} \left[ \langle (1+2K_h \eta), v \rangle + \bar{v} \langle \eta, (\partial_h K_h) \eta \rangle \right] = 0,
$$
where we have used $T_h \partial_u 1 = 0$, $\mathcal{M}_h^* 1 = (1+2K_h \eta)$, and integration by parts. Hence, $\delta M$ in (\ref{lin-mass}) is constant in time $t$. Due to the mass constraint $M(\eta) = 0$, it follows that $\delta M = 0$. 

Projection of (\ref{lin-2}) to $1$ yields
$$
\frac{d}{dt}  \langle (1+K_h \eta), w \rangle = -\langle (1+2K_h \eta), v \rangle + \bar{v} \langle 1, \partial_h S_h(\eta) \rangle  = -\delta M = 0,
$$
where we have used $\mathcal{M}_h 1 = 1 + K_h \eta$, $\mathcal{L}_h 1 = -(1+2K_h \eta)$, and $\langle 1, \partial_h S_h(\eta) \rangle = - \langle \eta, (\partial_h K_h) \eta \rangle$. Hence, $\delta Q$ in (\ref{lin-momentum-vert}) is constant in time $t$.

Projection of (\ref{lin-2}) to $\eta'$ yields
$$
\frac{d}{dt} \left[ \langle \eta', w \rangle - 2c \langle K_h \eta, v \rangle - 
c \bar{v} \langle \eta, (\partial_h K_h) \eta \rangle \right] = 0,
$$
where we have used $\mathcal{L}_h \eta' = 0$, $\mathcal{M}_h \eta' = \eta'$, and integration by parts. We have also used $\langle \eta', \partial_h S_h(\eta) \rangle = 0$ and $\langle \eta', K_h \eta \rangle = 0$ because $T_h^{-1}$ is skew-adjoint in $L^2_{\rm per}$. Hence, $\delta P$ in (\ref{lin-momentum}) is constant in time $t$.

Finally, we check the energy conservation. Expanding the energy $H(\zeta,\eta)$ in (\ref{energy-H-again}) gives the linear correction $\delta H$ to $H(\zeta,\eta)$ in the form:
\begin{align*}
\delta H &= 2 \langle c^2 K_h \eta + \eta (1+K_h \eta) + \frac{1}{2} K_h \eta^2, v \rangle + \bar{v} \left[ c^2 \langle \eta, (\partial_h K_h) \eta \rangle + \langle \eta^2, (\partial_h K_h) \eta \rangle \right] - 2 c \langle w, \eta' \rangle \\
&= 4 c^2 \langle K_h \eta, v \rangle - 2 c \langle w, \eta' \rangle + \bar{v} \left[ c^2 \langle \eta, (\partial_h K_h) \eta \rangle + \langle \eta^2, (\partial_h K_h) \eta \rangle + 2 B \langle \eta, (\partial_h K_h) \eta \rangle + 4 \pi B \right] \\
&= 2c \delta P,
\end{align*}
where we have used Babenko's equation (\ref{eq:Babenko}), the integral constraint (\ref{eq:DepthIdentity}), the definition (\ref{mean-term}), and the first-order corrections (\ref{lin-mass}) and (\ref{lin-momentum}). Since $\delta P$ is constant in time, $\delta H = 2c \delta P$ is also constant in time $t$.
\end{proof}

Proposition \ref{prop-lin-energy} shows that the linearized equations of motion (\ref{lin-1}) and (\ref{lin-2}) inherit conservation of $M(\eta)$, $Q(\zeta,\eta)$, $P(\zeta,\eta)$, and $H(\zeta,\eta)$ of the nonlinear equations (\ref{time-3}) and (\ref{time-4}). Since the constraint (\ref{time-constraint-expanded}) is equivalent to the energy conservation by Proposition \ref{prop-energy}, the linearized constraint (\ref{lin-constraint}) is equivalent to the time conservation of the second variation of the energy $H(\zeta,\eta)$ at the traveling wave solution with $\zeta = 0$ and the periodic profile $\eta$.
	
\subsection{Spectral problem}

Separating variables as 
$$
v(u,t) = e^{\lambda t} v(u), \quad w(u,t) = e^{\lambda t} w(u), \quad 
\bar{v}(t) = e^{\lambda t} \bar{v}, \quad (\delta B) = e^{\lambda t} b,
$$ 
we rewrite the linearized system of equations (\ref{lin-1}) and (\ref{lin-2}) in the form:
\begin{align}
-T_h \partial_u w &= \lambda \left( \mathcal{M}_h v - \bar{v} \eta' (\partial_h T_h^{-1}) \eta \right),  \label{spec-1} \\
\mathcal{L}_h v + \bar{v} \partial_h S_h(\eta) - 2b K_h \eta &=  \lambda \left( \mathcal{M}_h^* w + 2c T_h^{-1} v + c \bar{v} (\partial_h T_h^{-1}) \eta \right). & \label{spec-2}
\end{align}
The linearized constraint (\ref{lin-constraint}) 
is rewritten in the form
\begin{align}
2\pi b \left[ 1 + \frac{1}{2\pi} \langle \eta, (\partial_h K_h) \eta \rangle \right] &= \langle v, \partial_h S_h(\eta) \rangle 
+ \frac{1}{2} \bar{v} \langle (c^2 - 2B)  \eta - \eta^2, (\partial_h^2 K_h)\, \eta \rangle \notag \\
& \quad 
+ \lambda \left( \langle  w, \eta' (\partial_h T_h^{-1})\,\eta \rangle 
+ c \langle v, (\partial_h T^{-1}_h) \eta \rangle \right). 
\label{spec-constraint}
\end{align}
The system of equations (\ref{spec-1}), (\ref{spec-2}), and (\ref{spec-constraint}) defines the spectral problem for co-periodic 
perturbations. Since $H^1_{\rm per} \times H^1_{\rm per}$ is compactly embedded into $L^2_{\rm per} \times L^2_{\rm per}$, 
the linear operators $\mathcal{L}_h$ and $T_h \partial_u$ are unbounded in $L^2_{\rm per}$ with the domains in 
$H^1_{\rm per}$, and the linear operators $\mathcal{M}_h$, $\mathcal{M}_h^*$, and $T_h^{-1}$ are bounded in $L^2_{\rm per}$, 
the spectral problem admits only isolated eigenvalues $\lambda \in \mathbb{C}$. 

The following proposition identifies the constraints on the eigenfunction $(v,w) \in H^1_{\rm per} \times H^1_{\rm per}$ for a nonzero 
isolated eigenvalue $\lambda \in \mathbb{C} \backslash \{0\}$.

\begin{proposition}
\label{prop-orth}
Let  $(v,w) \in H^1_{\rm per} \times H^1_{\rm per}$ be an eigenfunction of the spectral problem (\ref{spec-1}) and (\ref{spec-2}) for $\lambda \neq 0$. 
Then, $\delta M = \delta Q = \delta P = 0$ in (\ref{lin-mass}), (\ref{lin-momentum-vert}), and (\ref{lin-momentum}). 
\end{proposition}

\begin{proof}
Projection of (\ref{spec-1}) to $1$ yields for $\lambda \neq 0$, 
\begin{equation*}
\langle (1+2K_h \eta), v \rangle + \bar{v} \langle \eta, (\partial_h K_h) \eta \rangle = 0,
\end{equation*}
which is the same as $\delta M = 0$ in (\ref{lin-mass}). Projection of (\ref{spec-2}) to $1$ yields for $\lambda \neq 0$, 
\begin{equation*}
\langle (1+K_h \eta), w \rangle = 0,
\end{equation*}
which is the same as $\delta Q = 0$ in (\ref{lin-momentum-vert}). Projection of (\ref{spec-2}) to $\eta'$ yields for $\lambda \neq 0$, 
\begin{equation*}
2c \langle K_h \eta, v \rangle - \langle \eta', w \rangle
+ c \bar{v} \langle \eta, (\partial_h K_h) \eta \rangle = 0,
\end{equation*}
which is the same as $\delta P = 0$ in (\ref{lin-momentum}). 
\end{proof}

\subsection{Multiplicity of the zero eigenvalue}

The main result of this section is the following theorem, which defines a sharp criterion that the zero eigenvalue admits a higher algebraic multiplicity 
compared to the quadruple multiplicity supported by the symmetries of the spectral problem, see Proposition \ref{prop-orth}. The criterion 
coincides with that for the deep fluid with $h = \infty$, see \cite{DP25}. The criterion is satisfied at each critical point of the horizontal momentum, 
which coincides with each critical point of the energy by Corollary \ref{prop-momentum-1-par}, in agreement with \cite{Z-S86}.

\begin{theorem}
	\label{theorem-multiplicity}
Assume the non-degeneracy condition ${\rm Ker}(\mathcal{L}_h) = {\rm span}(\eta')$, the solvability condition (\ref{smooth-cont}), and 
the additional condition
	\begin{equation}
	\label{IFT-condition}
	1 + \frac{1}{2\pi} \langle \eta, \partial_h K_h \eta \rangle \neq 0.
	\end{equation}
Then, the geometric multiplicity of $\lambda = 0$ in the spectral problem (\ref{spec-1})--(\ref{spec-2}) closed by the integral constraint (\ref{spec-constraint}) is equal to $2$, whereas the algebraic 
multiplicity of $\lambda = 0$ is equal to $4$ if and only if 
$\frac{d}{dc} \mathcal{P}(c,h(c)) \neq 0$, where $\mathcal{P}(c,h) = c \langle K_h \eta, \eta \rangle$ is given by (\ref{hor-momentum}) and $\frac{d}{dc} \mathcal{P}(c,h(c))$ is defined by (\ref{derivative-momentum-1-par}).
\end{theorem}

\begin{proof}
	The proof consists of three consecutive computations. \\
	
	\underline{The geometric multiplicity of $\lambda = 0$ is equal to $2$.} Setting $\lambda = 0$ in the spectral problem (\ref{spec-1})--(\ref{spec-2}) yields two nontrivial solutions:
$$
\mbox{\rm either} \quad (v,w) = (0,1) \quad \mbox{\rm or} \quad 
(v,w) = (\mathfrak{v},0),
$$ 
where $\mathfrak{v} \in H^1_{\rm per}$ is a general solution of the linear equation 
\begin{equation}
\label{lin-homog-eq}
\mathcal{L}_h \mathfrak{v} + \bar{v} \partial_h S_h(\eta) - 2b K_h \eta = 0.
\end{equation}
Both $\bar{v}$ and $b$ are obtained from $\mathfrak{v}$ by $\bar{v} = \frac{1}{2\pi} \langle 1, \mathfrak{v}\rangle$ and 
\begin{align}
2\pi \left[ 1 + \frac{1}{2\pi} \langle \eta, (\partial_h K_h) \eta \rangle \right] b &= \langle \mathfrak{v}, \partial_h S_h(\eta) \rangle 
+ \frac{1}{2} \bar{v} \langle (c^2 - 2B) \eta - \eta^2, (\partial_h^2 K_h) \eta \rangle,
\label{spec-constraint-zero}
\end{align}
where (\ref{spec-constraint-zero}) follows from (\ref{spec-constraint}) for $\lambda = 0$. Due to the constraint 
(\ref{IFT-condition}), the value of $b$ is uniquely defined by (\ref{spec-constraint-zero}). Since $\partial_h S(\eta)$ and $K_h \eta$ are orthogonal to $\eta'$ due to the skew-adjointness of $T_h^{-1}$ in $L^2_{\rm per}$, the non-degeneracy condition  ${\rm Ker}(\mathcal{L}_h) = {\rm span}(\eta')$ implies that $\partial_h S(\eta)$ and $K_h \eta$  belong to the range of $\mathcal{L}_h$, 
so that there exists an inhomogeneous solution $\mathfrak{v} \in H^1_{\rm per}$ of (\ref{lin-homog-eq}) for any $\bar{v}, b \in \mathbb{R}$. It follows from (\ref{derivative-1}), (\ref{derivative-2}), and (\ref{derivative-3}) that a general solution of (\ref{lin-homog-eq}) is given by 
\begin{equation}
\label{kernel}
\mathfrak{v} = \alpha_1 \eta' + \alpha_2 \partial_h \eta + \alpha_3 \partial_c \eta, 
\end{equation}
where $\alpha_1 \in \mathbb{R}$ is arbitrary and 
$\alpha_2,\alpha_3 \in \mathbb{R}$ are related to $\bar{v}$, $b$ by 
\begin{equation}
\label{kernel-alpha}
\alpha_2 = \bar{v}, \quad 
\alpha_3 = \frac{b - \bar{v} \partial_h B}{\partial_c B - c}.
\end{equation}
Proposition \ref{prop-2-parameters} under the condition (\ref{smooth-cont}) implies that $\partial_c B \neq c$ so that $\alpha_3$ is well-defined. 
We show that $\alpha_2 = \alpha_3 = 0$, which implies that the geometric multiplicity of $\lambda = 0$ is equal to $2$. 
By taking the mean value of (\ref{kernel}), we obtain 
$$
\bar{v} = \alpha_2 \overline{\partial_h \eta} + \alpha_3 \overline{\partial_c \eta}. 
$$
We use (\ref{derivative-constraints}) and obtain from (\ref{kernel-alpha}) that 
\begin{equation}
\label{technical-v}
\bar{v} = \frac{b - \bar{v} \partial_h B}{\partial_c B - c} h'(c).
\end{equation}
Substituting (\ref{kernel}) into (\ref{spec-constraint-zero}) 
yields 
\begin{align*}
2\pi \left[ 1 + \frac{1}{2\pi} \langle \eta, (\partial_h K_h) \eta \rangle \right] b &= \alpha_2 \langle \partial_h \eta, \partial_h S_h(\eta) \rangle 
+ \alpha_3 \langle \partial_c \eta, \partial_h S_h(\eta) \rangle \\
& \quad 
+ \frac{1}{2} \bar{v} \langle (c^2 - 2B) \eta - \eta^2, (\partial_h^2 K_h) \eta \rangle.
\end{align*}
Using (\ref{derivative-Ber-1}) and (\ref{kernel-alpha}), we obtain 
\begin{align*}
2 \pi \left[ 1 + \frac{1}{2\pi} \langle \eta, (\partial_h K_h) \eta \rangle \right] (b - \bar{v} \partial_h B)
&= \alpha_3 \langle \partial_c \eta, \partial_h S_h(\eta) \rangle,
\end{align*}
or equivalently, 
\begin{align*}
2\pi \left[ 1 + \frac{1}{2\pi} \langle \eta, (\partial_h K_h) \eta \rangle \right] \alpha_3 (\partial_c B - c) = \alpha_3 \langle \partial_c \eta, \partial_h S_h(\eta) \rangle.
\end{align*}
Eliminating $(\partial_c B)$ by using (\ref{derivative-Ber-2}) yields $-2\pi c \alpha_3 = 0$, which implies that $\alpha_3 = 0$. Hence, $b = \bar{v}\partial_h B$, which yields $\bar{v} = 0$ from (\ref{technical-v}), which is equivalent to $\alpha_2 = 0$.\\

\underline{The algebraic multiplicity of  $\lambda = 0$ is at least $4$.} 
We differentiate (\ref{spec-1})--(\ref{spec-2}) in $\lambda$ 
at $\lambda = 0$ and substitute
\begin{equation}
\label{kernel-alg}
v = \alpha \eta', \quad w = \beta, \quad \bar{v} = 0, \quad b = 0,
\end{equation}
with arbitrary $\alpha$ and $\beta$. This yields the linear inhomogeneous system of equations:
\begin{align}
-T_h \partial_u w_1 &=  \alpha \mathcal{M}_h \eta',  \label{gener-1} \\
\mathcal{L}_h v_1 + \overline{v}_1 \partial_h S_h(\eta) - 2 b_1 K_h \eta &=  \beta \mathcal{M}_h^* 1 + 2c \alpha T_h^{-1} \eta', & \label{gener-2}
\end{align}
where $\overline{v}_1 = \frac{1}{2\pi} \langle 1, v_1 \rangle$ and $b_1$ 
are derivatives of $\overline{v}$ and $b$ in $\lambda$ at $\lambda = 0$. The system (\ref{gener-1})--(\ref{gener-2}) is closed by the constraint, which is obtained from differentiating (\ref{spec-constraint}) in $\lambda$ at $\lambda = 0$:
\begin{align}
2\pi b_1 \left[ 1 + \frac{1}{2\pi} \langle \eta, (\partial_h K_h) \eta \rangle \right] &= \langle v_1, \partial_h S_h(\eta) \rangle 
+ \frac{1}{2} \overline{v}_1 \langle (c^2 - 2B)  \eta - \eta^2, (\partial_h^2 K_h)\, \eta \rangle \notag \\
& \quad 
+ (\beta + c \alpha) \langle \eta', (\partial_h T^{-1}_h) \eta \rangle,
\label{spec-constraint-alg}
\end{align}
where the last projection is generally nonzero. The solvability condition for  (\ref{gener-1}) is obtained by its projection to ${\rm Ker}(T_h \partial_u) = {\rm span}(1)$. The solvability condition is satisfied for arbitrary $\alpha$  and $\beta$ since 
\begin{align*}
\alpha \langle 1, \mathcal{M}_h \eta' \rangle = \alpha \langle (1+2K_h \eta), \eta' \rangle = 0.
\end{align*}
The solvability condition for (\ref{gener-2}) is obtained by its projection to ${\rm Ker}(\mathcal{L}_h) = {\rm span}(\eta')$. The solvability condition  is satisfied for arbitrary $\alpha$ and $\beta$ since
\begin{align*}
\langle \eta', \beta \mathcal{M}_h^* 1 + 2c \alpha T_h^{-1} \eta' \rangle 
= \beta \langle \eta', 1 \rangle + 2 c \alpha \langle \eta', K_h \eta \rangle = 0. 
\end{align*}
Hence, for every $(\alpha,\beta)$, we can construct the generalized eigenvectors $(v_1,w_1) \in H^1_{\rm per} \times H^1_{\rm per}$ from solutions of the linear inhomogeneous system (\ref{gener-1})--(\ref{gener-2}). By using (\ref{derivative-2}), (\ref{derivative-3}), $\mathcal{M}_h \eta' = \eta'$, $\mathcal{M}_h^* 1 = 1 + 2 K_h \eta$, and $\mathcal{L}_h 1 = -(1+2K_h \eta)$, we can obtain 
the exact solution of (\ref{gener-1})--(\ref{gener-2})  in the form:
\begin{align}
\label{sol-v-1}
v_1 = \bar{v}_1 \partial_h \eta + \frac{b_1 - \bar{v}_1 \partial_h B + c \alpha}{\partial_c B - c} \partial_c \eta - \beta, \quad w_1 = -\alpha T_h^{-1} \eta.
\end{align}
By taking the mean value and using (\ref{derivative-constraints}), we obtain 
similarly to (\ref{technical-v}), 
\begin{align}
\label{bar-v-1}
\bar{v}_1 = \frac{b_1 - \bar{v}_1 \partial_h B + c \alpha}{\partial_c B - c} h'(c) - \frac{\beta}{1 - \overline{\partial_h \eta}}.
\end{align}
By using (\ref{spec-constraint-alg}), we obtain 
\begin{align*}
2\pi b_1 \left[ 1 + \frac{1}{2\pi} \langle \eta, (\partial_h K_h) \eta \rangle \right] &= \frac{b_1 - \bar{v}_1 \partial_h B + c \alpha}{\partial_c B - c} \langle \partial_c v, \partial_h S_h(\eta) \rangle 
- c \alpha \langle \eta, (\partial_h K_h) \eta \rangle \\
& \quad + \bar{v}_1 \langle \partial_h \eta, \partial_h S_h(\eta) \rangle 
+ \frac{1}{2} \overline{v}_1 \langle (c^2 - 2B)  \eta - \eta^2, (\partial_h^2 K_h)\, \eta \rangle,
\end{align*}
where the $\beta$ term is canceled because 
$$
\langle 1, \partial_h S_h(\eta) \rangle = -\langle \eta, \partial_h K_h \eta \rangle \quad \mbox{\rm and} \quad 
\langle \eta', (\partial_h T^{-1}_h) \eta \rangle = -\langle \eta, \partial_h K_h \eta \rangle.
$$
By using (\ref{derivative-Ber-1}) and (\ref{derivative-Ber-2}), we obtain after simplification,
\begin{align*}
2\pi \left( b_1 - \bar{v}_1 \partial_h B  - \frac{b_1 - \bar{v}_1 \partial_h B + c \alpha}{\partial_c B - c} \partial_c B \right) = 0,
\end{align*}
which yields $b_1 = \bar{v}_1 \partial_h B - \alpha \partial_c B$. By using 
(\ref{bar-v-1}), we derive 
\begin{align}
\label{kernel-alg-2-param}
\bar{v}_1 = -\alpha h'(c) - \frac{\beta}{1 - \overline{\partial_h \eta}} \quad \mbox{\rm and} \quad 
b_1 = -\alpha (\partial_c B + h'(c) \partial_h B) - \frac{\beta \partial_h B}{1 - \overline{\partial_h \eta}}.
\end{align}
Substituting these expressions into (\ref{sol-v-1}) yields the exact solution to (\ref{gener-1})--(\ref{gener-2}) for arbitrary $\alpha$ and $\beta$:
\begin{align}
\label{kernel-alg-2}
v_1 = - \alpha \left( \partial_c \eta + h'(c) \partial_h \eta \right) 
- \beta \left( 1 + \frac{\partial_h \eta}{1 - \overline{\partial_h \eta}} \right), \quad w_1 = -\alpha T_h^{-1} \eta.
\end{align}
Thus, we obtain two generalized eigenvectors in (\ref{kernel-alg-2}) associated with 
the two eigenvectors in (\ref{kernel-alg}), suggesting the Jordan block structure of at least two double blocks for $\lambda = 0$.\\

\underline{Threshold for the quadruple algebraic multiplicity of  $\lambda = 0$.} 
We continue the Jordan chain of generalized eigenvectors beyond 
(\ref{kernel-alg}) and (\ref{kernel-alg-2}). The second generalized eigenvectors are defined from the linear inhomogeneous system of equations:
	\begin{align}
	-T_h \partial_u w_2 &= \mathcal{M}_h v_1 - \bar{v}_1 \eta' (\partial_h T_h^{-1}) \eta,  \label{second-1} \\
	\mathcal{L}_h v_2 + \overline{v}_2 \partial_h S_h(\eta) - 2 b_2 K_h \eta &=  \mathcal{M}_h^* w_1 + 2c T_h^{-1} v_1 + c \bar{v}_1 (\partial_h T_h^{-1}) \eta, & \label{second-2}
	\end{align}
where $\overline{v}_2 = \frac{1}{2\pi} \langle 1, v_2 \rangle$ and $b_2$ 
are equal to the second derivatives of $\overline{v}$ and $b$ in $\lambda$ at $\lambda = 0$ divided by $2$.
The system (\ref{second-1})--(\ref{second-2}) is closed by the constraint
	\begin{align}
	2\pi b_2 \left[ 1 + \frac{1}{2\pi} \langle \eta, (\partial_h K_h) \eta \rangle \right] &= \langle v_2, \partial_h S_h(\eta) \rangle 
	+ \frac{1}{2} \overline{v}_2 \langle (c^2 - 2B)  \eta - \eta^2, (\partial_h^2 K_h)\, \eta \rangle,
	\label{spec-second}
	\end{align}
	where we have used 
	$$
	\langle w_1, \eta' (\partial_h T_h^{-1}) \eta \rangle = \langle v_1, (\partial_h T_h^{-1}) \eta \rangle = 0.
	$$
Since projection of (\ref{gener-2}) to $1$ yields
	\begin{align*}
	-\langle (1 + 2K_h \eta), v_1 \rangle - \bar{v}_1 \langle \eta, (\partial_h K_h) \eta \rangle = 2 \pi \beta, 
	\end{align*}
	which can be checked directly from (\ref{kernel-alg-2-param}) and (\ref{kernel-alg-2}), we obtain the solvability condition for (\ref{second-1}) by its projection to ${\rm Ker}(T_h \partial_u) = {\rm span}(1)$:
\begin{align*}
0 = \langle 1, \mathcal{M}_h v_1 - \bar{v}_1 \eta' (\partial_h T_h^{-1}) \eta \rangle = \langle (1+2K_h \eta), v_1 \rangle 
+ \bar{v}_1 \langle \eta, (\partial_h K_h) \eta \rangle = -2\pi \beta.
\end{align*}	
The solvability condition is satisfied if and only if $\beta = 0$, which means that the double 
Jordan block associated with the $\beta$ terms cannot be continued beyond the first generalized eigenvector. Next, we obtain the solvability condition for (\ref{second-2}) by its projection to ${\rm Ker}(\mathcal{L}_h) = {\rm span}(\eta')$:
	\begin{align*}
	0 &= \langle \eta', \mathcal{M}_h^* w_1 + 2c T_h^{-1} v_1 + c \bar{v}_1 (\partial_h T_h^{-1}) \eta \rangle.
	\end{align*}	
Substituting $\bar{v}_1$ from (\ref{kernel-alg-2-param}), 
and $v_1,w_1$ from (\ref{kernel-alg-2}) with $\beta = 0$ yields
	\begin{align*}	
0 &= \alpha \left[	\langle \eta', -\mathcal{M}_h^* T_h^{-1} \eta - 2c T_h^{-1} (\partial_c \eta + h'(c) \partial_h \eta) - c h'(c) (\partial_h T_h^{-1}) \eta \rangle \right] \\
&= \alpha \left[ \langle K_h \eta, \eta \rangle - -c \langle 1, \partial_c \eta \rangle - c h'(c) \langle 1, \partial_h \eta \rangle \right] \\
	&= \alpha \left[ \langle K_h \eta, \eta \rangle - 2 \pi c h'(c) \right] \\
	&= \alpha \frac{d}{dc} \mathcal{P}(c,h(c)),
	\end{align*}
	where we have used (\ref{derivative-momentum}) and (\ref{derivative-momentum-1-par}). If $\frac{d}{dc} \mathcal{P}(c,h(c)) \neq 0$, the solvability condition is satisfied if and only if $\alpha = 0$, which means that the double Jordan block associated with the $\alpha$ terms cannot be continued beyond the first generalized eigenvector. If  $\frac{d}{dc} \mathcal{P}(c,h(c)) = 0$, the solvability condition is satisfied and 
the second generalized eigenvector exists in the system (\ref{second-1}), (\ref{second-2}), and (\ref{spec-second}) with $\alpha \neq 0$.
\end{proof}

\begin{corollary}
	\label{cor-multiplicity}
	Assume that the periodic profile $\eta \in C^{\infty}_{\rm per}$ is even in $u$ on $\mathbb{T}$. 	Under the conditions of Theorem \ref{theorem-multiplicity}, the algebraic multiplicity of $\lambda = 0$ is equal to at least $6$ if $\frac{d}{dc} \mathcal{P}(c,h(c)) = 0$.
\end{corollary}

\begin{proof}
	 The third generalized eigenvectors are defined from the linear inhomogeneous system of equations:
	\begin{align}
	-T_h \partial_u w_3 &= \mathcal{M}_h v_2 - \bar{v}_2 \eta' (\partial_h T_h^{-1}) \eta,  \label{third-1} \\
	\mathcal{L}_h v_3 + \overline{v}_3 \partial_h S_h(\eta) - 2 b_3 K_h \eta &=  \mathcal{M}_h^* w_2 + 2c T_h^{-1} v_2 + c \bar{v}_2 (\partial_h T_h^{-1}) \eta, & \label{third-2}
	\end{align}
where $\overline{v}_3 = \frac{1}{2\pi} \langle 1, v_3 \rangle$ and $b_3$ 
are equal to the third derivatives of $\overline{v}$ and $b$ in $\lambda$ at $\lambda = 0$ divided by $6$.
Projection of (\ref{second-2}) to $1$ yields the relation 
\begin{equation*}
\langle (1+2K_h \eta), v_2 \rangle +  \bar{v}_2 \langle \eta, (\partial_h K_h) \eta \rangle = 0,
\end{equation*}
from which the solvability condition for (\ref{third-1}) is satisfied:
\begin{equation*}
\langle 1, \mathcal{M}_h v_2 - \bar{v}_2 \eta' (\partial_h T_h^{-1}) \eta\rangle = 0.
\end{equation*}
In addition,  $v_2$ is odd and $w_2$ is even in $u$ if $\eta$ is even in $u$. Hence, the solvability condition for (\ref{third-2}) is satisfied,
\begin{equation*}
\langle \eta', \mathcal{M}_h^* w_2 + 2c T_h^{-1} v_2 + c \bar{v}_2 (\partial_h T_h^{-1}) \eta \rangle = 0,
\end{equation*}
so that the third generalized eigenvector exists in the system (\ref{third-1}), (\ref{third-2}), and the associated constraint. 
Taking into account the double Jordan block for the $\beta$ terms and 
the (at least) quadruple Jordan block for the $\alpha$ terms, this yields 
the algebraic multiplicity of $\lambda = 0$ of at least $6$.
\end{proof}

\section{Small-amplitude expansion of the wave profile}
\label{sec-small}

Here, we compare the small-amplitude expansion of the traveling wave profile obtained via the conformal transformation with the recent computations 
based on other formulations of Euler's equations. To enable comparison with \cite{BMV23,CDT22}, 
we augment Babenko's equation (\ref{eq:Babenko}) and the integral constraint (\ref{eq:DepthIdentity}) with the nonzero mass invariant $M(\eta)$. This yields the system of equations for the traveling wave profile $\eta \in C^{\infty}_{\rm per}$:
\begin{equation}
\label{eq:Babenko-small}
(c^{2} - 2B) K_h \eta - \eta
= \frac12 K_h \eta^{2} + \eta\,K_h\eta - \frac{1}{2\pi} M(\eta)
\end{equation}
and
\begin{equation}
\label{eq:DepthIdentity-small}
B + \frac{1}{2\pi} M(\eta) = \frac{1}{4\pi} \langle  (c^{2} - 2B) \eta - \eta^2, (\partial_h K_h) \eta \rangle.
\end{equation}
Since the steady traveling waves are symmetric with respect to the crest,  the wave profile $\eta \in C^{\infty}_{\rm per}$ is considered to be an \emph{even} function of $u$.  The Fourier expansion of the even wave profile 
contains only cosine harmonics:
\begin{equation}\label{eq:EtaCosineForm}
\eta(u;h,a)
= \sum_{n = 0}^{\infty} \eta^{(n)}(h,a) \cos(nu),
\end{equation}
where $a > 0$ is the (small) amplitude parameter. 
In the small–amplitude expansion, we define
\begin{equation}
\label{eq:EtaExpansion}
\eta(u;h,a)
= a\cos(u)
+ a^{2} \left[ A_0(h) + A_2(h)\cos(2u) \right] 
+ a^{3}A_3(h)\cos(3u)
+ \mathcal{O}(a^{4}),
\end{equation}
with the wave speed $c$ expanded as 
\begin{equation}\label{eq:c2Expansion}
c^{2}(h,a)
= c_{0}^{2}(h)
+ a^{2} c_{2}(h)
+ O(a^{4}),
\end{equation}
the Bernoulli constant $B$ expanded as 
\begin{equation}
\label{eq:BernouliExpansion}
B = a^2 B_2(h) + \mathcal{O}(a^4),
\end{equation}
and the mass invariant $M(\eta)$ expanded as 
\begin{equation}
\label{mass}
M(\eta) = a^2 \left( \langle 1, \eta_2 \rangle + \langle \eta_1, K_h\eta_1 \rangle \right) + \mathcal{O}(a^4).
\end{equation}
Comparison between (\ref{eq:EtaCosineForm}) and (\ref{eq:EtaExpansion}) yields 
\begin{align*}
\left\{ \begin{array}{l} 
\eta^{(0)}(h,a) = a^2 A_0(h) + \mathcal{O}(a^4),\\
\eta^{(1)}(h,a) = a, \\
\eta^{(2)}(h,a) = a^2 A_2(h) + \mathcal{O}(a^4), \\ 
\eta^{(3)}(h,a) = a^3 A_3(h) + \mathcal{O}(a^5),
\end{array} \right.
\end{align*}
where $\eta^{(1)}(h,a) = a$ is not expanded in powers of $a$ due to an arbitrary normalization of the amplitude parameter $a$. We also note that even powers of $a$ enter even Fourier cosine modes and odd powers of $a$ enter odd Fourier cosine modes.
	
For a fixed conformal depth $h>0$, the operator $K_h$ is defined by its action on Fourier modes \eqref{eq:KhDefinition}. Applied to the cosine modes, it gives
	\begin{equation}\label{eq:KhDef}
	K_h \cos(nu) = n \coth(nh)\, \cos(nu), \qquad n\in\mathbb{N}.
	\end{equation}
Since the eigenfunctions $\{\cos(nu)\}_{n\ge 1}$ form an orthogonal basis for the space of even $2\pi$–periodic functions, the operator $K_h$ is diagonalized by this Fourier basis. Differentiating \eqref{eq:KhDef} with respect to $h$ yields
	\begin{equation}\label{eq:dKh}
	(\partial_h K_h)\cos(nu)
	= -\, \frac{n^{2}}{\sinh^{2}(nh)} \cos(nu),
	\qquad n\in\mathbb{N}.
	\end{equation}
Thus, $(\partial_h K_h)$ is again diagonal in the Fourier cosine basis, with real negative eigenvalues. For further use, we record the standard trigonometric identities
	\begin{align}
	\begin{cases} 
	2 \cos(u) \cos(u) =  1 + \cos 2u, \\
	2 \cos(u) \cos(2u) =  \cos u + \cos 3 u, 
	\end{cases}
	\label{cos-3}
	\end{align}
which we use in the computations of the small-amplitude expansion. 

In what follows, we compute the small-amplitude (Stokes) expansions and compare them with the outcomes of computations in 
\cite{BMV23}, \cite{CDT22}, and \cite{Zakharov96}. Justification 
of the Stokes expansion (\ref{eq:EtaExpansion}) is a classical result and hence 
it will be omitted.

\subsection{Computation of Stokes expansion from (\ref{eq:Babenko-small}) and (\ref{eq:DepthIdentity-small})}
	
We write formally $\eta = a \eta_1 + a^2 \eta_2 + a^3 \eta_3 + \mathcal{O}(a^4)$ and proceed with computations. 

$\phantom{t}$\\	
\underline{Order $\mathcal{O}(a)$:} Equation (\ref{eq:Babenko-small}) reduces to a linear eigenvalue problem for the bifurcating mode $\eta_1(u) = \cos(u)$:
	\begin{equation}
	\label{inhom-0}
	(c^2_0 K_h - 1) \eta_1 = 0.
	\end{equation} 
	Using the diagonal action of $K_h$ on cosine modes, given by (\ref{eq:KhDef}) for $n = 1$, we obtain 
	\[
	c_0^2(h)\,\coth(h) = 1,
	\]
	which yields the leading–order dispersion relation
	\begin{equation}\label{eq:c0Dispersion}
	c_0^2(h) = \tanh(h).
	\end{equation}
	
$\phantom{t}$\\	
\underline{Order $\mathcal{O}(a^2)$:} Equation (\ref{eq:Babenko-small}) reduces to the linear inhomogeneous equation 
	\begin{equation}
	\label{inhom-1}
	(c_0^2 K_h - 1) \eta_2 + \frac{1}{2\pi} \langle 1, \eta_2 \rangle = \frac{1}{2} K_h \eta_1^2 + \eta_1 K_h \eta_1 - \frac{1}{2\pi} \langle \eta_1, K_h\eta_1 \rangle.
	\end{equation}
By using (\ref{eq:KhDef}) and the first line in (\ref{cos-3}), we obtain
	$$
	\frac{1}{2} K_h \eta_1^2 + \eta_1 K_h \eta_1 = \frac{1}{2} \coth(2h) \cos(2u) + \frac{1}{2} \coth(h) (1 + \cos(2u)).
	$$
Since
	$$
	\frac{1}{2\pi} \langle \eta_1, K_h\eta_1 \rangle = \frac{1}{2} \coth(h),
	$$
the right-hand side of (\ref{inhom-1}) contains only projection to $\cos(2u)$:
$$
	\frac{1}{2} K_h \eta_1^2 + \eta_1 K_h \eta_1 - \frac{1}{2\pi} \langle \eta_1, K_h\eta_1 \rangle = \frac{1}{2} [\coth(2h) + \coth(h)] \cos(2u).
	$$
	Solving (\ref{inhom-1}) with $\eta_2(u;h) = A_0(h) + A_2(h) \cos(2u)$ yields 
	an arbitrary $A_0(h)$ and a uniquely defined $A_2(h)$ by 
	\begin{equation}\label{eq:A2h}
	A_2(h) = \frac{\coth(2h) + \coth(h)}{2 (2 \tanh(h) \coth(2h) - 1)}
	= \frac{c_0^4(h) + 3}{4\,c_0^6(h)},
	\end{equation}
	where we have transformed (\ref{eq:A2h}) by using the identities:
	\[
	\coth(h) = \frac{1}{c_0^2(h)}, \quad \coth(2h) = \frac{1 + c_0^4(h)}{2c_0^2(h)}.
	\]
The expansion of the mass $M(\eta)$ in (\ref{mass}) is given by 
	\begin{align}
	\label{mean-value}
	\frac{1}{2\pi} M(\eta) = a^2 \left[ A_0(h) + \frac{1}{2} \coth(h) \right] + \mathcal{O}(a^4).
	\end{align}	
Although $A_0(h)$ is not defined by (\ref{inhom-1}), we get a relation between $A_0(h)$ and $B_2(h)$ from the integral constraint (\ref{eq:DepthIdentity-small}) at the order of $\mathcal{O}(a^2)$. 
By using (\ref{eq:dKh}), we obtain
\[
(\partial_h K_h)\eta
=
-a\,\csch^2(h) \cos(u)
-4a^2 A_2 \csch^2(2h)\cos(2u)
+ \mathcal{O}(a^3), 
\]
from which we derive
\begin{align}
\label{exp-1}
\frac{c^2}{2} \langle \eta, (\partial_h K_h)\eta \rangle  = -\frac{\pi}{2} c_0^2  \csch^2(h) a^2 + \mathcal{O}(a^4).
\end{align}
Collecting terms of the order of $\mathcal{O}(a^2)$ in (\ref{eq:DepthIdentity-small}), 
(\ref{eq:BernouliExpansion}), and (\ref{mean-value}) yields the relation between $A_0(h)$ and $B_2(h)$:
\begin{equation}
\label{def-A0}
A_0(h) + B_2(h) = -\frac{1}{2} \coth(h) -\frac{1}{2 \sinh(2h)},
\end{equation}
where we have used (\ref{eq:c0Dispersion}). 
	
$\phantom{t}$\\	
\underline{Order $\mathcal{O}(a^3)$:} Equation (\ref{eq:Babenko-small}) reduces to 
	the linear inhomogeneous equation 
	\begin{equation}
	\label{inhom-2}
(c_0^2 K_h - 1) \eta_3 + c_2 K_h \eta_1 = K_h \eta_1 \eta_2 + \eta_1 K_h \eta_2 + \eta_2 K_h \eta_1 + 2 B_2 K_h \eta_1,
	\end{equation}
	where we have used expansions (\ref{eq:c2Expansion}), (\ref{eq:BernouliExpansion}), and (\ref{mass}). The right–hand side of (\ref{inhom-2}) contains only projections to $\cos u$ and $\cos 3u$, which have zero mean. By using (\ref{eq:KhDef}) and the second line in (\ref{cos-3}), we obtain 
	\begin{align*}
	\begin{cases}
	K_h \eta_1\eta_2 = \frac{1}{2} (2 A_0+A_2) \coth(h) \cos(u)
	+ \frac{3}{2} A_2 \coth(3h) \cos(3u), \\
	\eta_2 K_h\eta_1 = \frac{1}{2} (2 A_0+A_2) \coth(h) \cos(u)
	+ \frac{1}{2} A_2 \coth(h) \cos(3u), \\
	\eta_1 K_h\eta_2 = A_2 \coth(2h) \cos(u) + A_2 \coth(2h) \cos(3u), 
	\end{cases}
	\end{align*}
	which yields
	\begin{align*}
K_h \eta_1\eta_2 +\eta_2 K_h\eta_1+\eta_1 K_h\eta_2 &= 
[ (2 A_0+A_2) \coth(h) + A_2 \coth(2h) ] \cos(u) \\
& \quad + \frac{1}{2} A_2 [\coth(h) + 2 \coth(2h) + 3 \coth(3h)] \cos(3u).
	\end{align*}
	The Fourier mode $\cos u$ is a solution of the homogeneous equation 
	(\ref{inhom-0}). To eliminate non-periodic solutions of (\ref{inhom-2}), 
	we remove it by the choice of $c_2$:
	\begin{equation*}
	c_2(h)
	= 2 B_2(h) + 2A_0(h) + A_2(h)
	+ A_2(h)\,\frac{\coth(2h)}{\coth h}.
	\end{equation*}
By using  (\ref{eq:A2h}) and (\ref{def-A0}), we find the explicit expression for $c_2(h)$:
\begin{align}
c_2(h) &= - \frac{1}{\sinh(2h)} -\frac{1}{\tanh(h)} + \frac{(3+\tanh^2(h))^2}{8 \tanh^3(h)} \notag \\
&= \frac{9 - 6\tanh^2(h) + 5 \tanh^4(h)}{8\tanh^3(h)}.
\label{eq:c2Condition}
\end{align}
	Solving (\ref{inhom-2}) with $\eta_3(u;h) = A_3(h) \cos(3u)$ by using (\ref{eq:KhDef}) and (\ref{eq:c0Dispersion}), we obtain  
	$$
	A_3(h) = \frac{1}{2} A_2(h) \frac{\coth(h) + 2 \coth(2h) + 3 \coth(3h)}{
		3\tanh(h) \coth(3h) - 1}
	$$
	Using
	\[
	\coth(3h) = \frac{\coth(h) \coth(2h) + 1}{\coth(2h) + \coth(h)} = \frac{1+3c_0^4(h)}{c_0^2(h) (3+c_0^4(h))},
	\]
	we get 
	\[
	\coth(h) + 2 \coth(2h) + 3 \coth(3h)  = \frac{c_0^8(h) + 14 c_0^4(h) + 9}{c_0^2(h) (3+c_0^4(h))}
	\]
	and 
	\[
	3\tanh(h) \coth(3h) - 1 =\frac{8c_0^4(h)}{3+c_0^4(h)}.
	\]
	which yields
	\begin{equation*}
	A_3(h)
	= \frac{\big(c_0^4(h) + 3\big)\big(c_0^8(h) + 14c_0^4(h) + 9\big)}
	{64\,c_0^{12}(h)}.
	\end{equation*}    
The integral constraint (\ref{eq:DepthIdentity-small}) contains no nonzero terms 
at the order of $\mathcal{O}(a^3)$ due to (\ref{eq:BernouliExpansion}), (\ref{mass}), (\ref{exp-1}), and 
\begin{align*}
-\frac12 \langle \eta^2, (\partial_h K_h) \eta  \rangle 
= \mathcal{O}(a^4). 
\end{align*}
This completes the construction of the leading-order expansions (\ref{eq:EtaExpansion}), (\ref{eq:c2Expansion}), (\ref{eq:BernouliExpansion}), and (\ref{mass}) up to the order of $\mathcal{O}(a^4)$ from Babenko's equation (\ref{eq:Babenko-small}) and the integral constraint  (\ref{eq:DepthIdentity-small}).

\subsection{Comparison with \cite{BMV23}}

The Stokes expansion is obtained in \cite{BMV23} by using the Craig--Sulem (CS) formulation of Euler's equations in physical coordinates. The traveling wave profiles are obtained for the zero value of the Bernoulli constant: $B = 0$. Hence, we set $B_2(h) = 0$ in (\ref{def-A0}) 
and obtain 
$$
A_0(h) = -\frac{1}{2} \coth(h) -\frac{1}{2 \sinh(2h)}.
$$
The mean value of the surface elevation in the physical coordinate 
$x$ is computed by the chain rule $dx = (1+K_h \eta) du$ from  (\ref{xi-expressions}) as $\eta^{[0]} := \frac{1}{2\pi} M(\eta)$. Due to (\ref{mean-value}), the second-order correction term $\eta^{[0]}_2$ recovers the expression in \cite[Eq. (2.6)]{BMV23}:
$$
\eta^{[0]}_2 = A_0(h) + \frac{1}{2} \coth(h) = -\frac{1}{2 \sinh(2h)} = \frac{c_0^4(h) - 1}{4 c_0^2(h)},
$$
where we have used $c_0^2(h) = \tanh(h)$. By using (\ref{xi-expressions}), we get 
$$
\xi - u = a \coth(h) \sin(u) + a^2 \coth(2h) \sin(2u) + \mathcal{O}(a^3),
$$
which is inverted as 
$$
u = x - a \coth(h) \sin(x) + \mathcal{O}(a^2). 
$$
The Stokes expansion (\ref{eq:EtaExpansion}) is rewritten as 
$$
\eta = a \cos(x) + a^2 \left( \eta^{[0]}_2 + \eta^{[2]}_2 \cos(2x) \right)+ \mathcal{O}(a^3),
$$
where 
$$
\eta^{[2]}_2 = A_2(h) - \frac{1}{2} \coth(h) = \frac{3 - c_0^4(h)}{4 c_0^6(h)},
$$
in agreement with \cite[Eq. (2.7)]{BMV23}. Using (\ref{eq:DepthConstraint}), we obtain 
\begin{equation}
\label{exp-h}
h = h_0 + a^2 A_0(h) + \mathcal{O}(a^4),
\end{equation}
in terms of the physical depth $h_0$. By using the expansion (\ref{eq:c2Expansion}) with the correction terms (\ref{eq:c0Dispersion}) and (\ref{eq:c2Condition}), we express the speed $c$ in terms of the physical depth as 
\begin{align*}
c^2 &= c_0^2(h) + a^2 c_2(h) + \mathcal{O}(a^4) \\
&= \tanh(h) + a^2 \frac{9 - 6 \tanh^2(h)+ 5 \tanh^4(h)}{8 \tanh^3(h)} + \mathcal{O}(a^4) \\
&= \tanh(h_0) + a^2 \left[ {\rm sech}^2(h_0) A_0(h_0) + \frac{9 - 6 \tanh^2(h_0) + 5 \tanh^4(h_0)}{8 \tanh^3(h_0)} \right] + \mathcal{O}(a^4),
\end{align*}
from which we recover the expression in \cite[Eq. (2.8)]{BMV23}:
\begin{align*}
c &= c_0(h_0) + \frac{a^2}{2 c_0(h_0)} \left[ 
(1 - c_0^4) \eta^{[0]}_2 - \frac{1}{2} \coth(h_0) {\rm sech}^2(h_0) + \frac{9 - 6 \tanh^2(h_0) + 5 \tanh^4(h_0)}{8 \tanh^3(h_0)} \right] + \mathcal{O}(a^4) \\
&= c_0(h_0) + \frac{a^2}{2 c_0(h_0)} \left[ 
(1 - c_0^4) \eta^{[0]}_2 + \frac{9 - 10 \tanh^2(h_0) + 9 \tanh^4(h_0)}{8 \tanh^3(h_0)} \right] + \mathcal{O}(a^4).
\end{align*}
Thus, the Stokes expansions (\ref{eq:EtaExpansion}), (\ref{eq:c2Expansion}), (\ref{eq:BernouliExpansion}), and (\ref{mass}) obtained via the conformal transformation recover the Stokes expansion from \cite{BMV23} obtained via the CS formulation of Euler's equations. 

\subsection{Comparison with \cite{CDT22}}

The Stokes expansion is obtained in \cite{CDT22} by using the 
Ablowitz--Fokas--Musslimani (AFM) formulation of Euler's equations in physical coordinates under an integral constraint. The traveling wave profiles are obtained for the zero value of the mass: $M(\eta) = 0$. 
This convention is used to normalize the physical depth $h_0$ by enforcing the mean value of the surface elevation in the physical coordinate $x$ to zero: $\eta^{[0]} = 0$. From  (\ref{mean-value}) and (\ref{def-A0}), this uniquely defines that 
\begin{equation}
\label{exp-B}
A_0(h) = -\frac{1}{2} \coth(h), \quad B_2(h) = -\frac{1}{2 \sinh(2h)}.
\end{equation}
The expansion (\ref{exp-h}) yields now 
$$
h = h_0 - \frac{1}{2} a^2 \coth(h) + \mathcal{O}(a^4),
$$
from which we derive
\begin{align}
\label{exp-c}
c = c_0(h_0) + \frac{a^2}{2 c_0(h_0)} \left[ 
\frac{9 - 10 \tanh^2(h_0) + 9 \tanh^4(h_0)}{8 \tanh^3(h_0)} \right] + \mathcal{O}(a^4),
\end{align}
which is different from the expression in \cite[Eq. (A.2b)]{CDT22} since 
\begin{align*}
\frac{6 + 2 \cosh(2h_0) + \cosh(4 h_0)}{8 \sinh^3(h_0) \cosh(h_0)}
&= \frac{9 + 12 \sinh^2(h_0) + 8\sinh^4(h_0)}{8 \tanh^3(h_0) \cosh^4(h_0)}  \\
&= \frac{9 - 6 \tanh^2(h_0) + 5\tanh^4(h_0)}{8 \tanh^3(h_0)} \\
&\neq \frac{9 - 10 \tanh^2(h_0) + 9 \tanh^4(h_0)}{8 \tanh^3(h_0)}.
\end{align*}

The discrepancy is explained as follows. It follows from \cite[Eqs. (A.1) and (A.2)]{CDT22} that the Bernoulli constant 
is also set to zero, $B = 0$, which is impossible in our normalization of the traveling wave profiles. However, it is clear from \cite[Eqs. (A.3b) and (A.4e)]{CDT22} that the speed $c$ is not defined by the expansion \cite[Eqs. (A.2)]{CDT22} because the velocity potential $\psi$ has a linear growth in $x$. This is also impossible in our normalization since the velocity potential $\psi$ is required to be a periodic function of $u$ (and hence, of $x$) in (\ref{Bernoulli}). 

The following proposition establishes the equivalence between the outcomes of the AFM method used in \cite{CDT22} and the outcomes of our method based on the conformal transformation.

\begin{proposition}
    \label{prop-comparison}
The speed in the AFM method used in \cite{CDT22}  coincides with the bifurcation parameter $\gamma$ used in Babenko's equation (\ref{BE}).
\end{proposition}

\begin{proof}
    We recall the Galilean transformation for Euler's equations (\ref{euler-1-2}). Let $\eta = \eta(x-\gamma t)$ be a $2\pi$-periodic traveling wave with speed $\gamma$ and $\varphi(x,y,t) = \alpha (x - \gamma t) + \tilde{\varphi}(x-\gamma t, y)$ be the velocity potential with 
    linear growth in $x$ given by $\alpha$ and the $2\pi$-periodic profile $\tilde{\varphi}(\cdot,y)$ with zero mean. 
    Euler's equations  (\ref{euler-1-2}) are rewritten in the form:
 \begin{align*}
 y = \eta(x,t) : \quad 
 \begin{cases}
-(\gamma - \alpha) \eta' + \tilde{\varphi}_x \eta' - \tilde{\varphi}_y = 0,  \\
B -(\gamma - \alpha) \tilde{\varphi}_x + \frac{1}{2} (\tilde{\varphi}_x)^2 + \frac{1}{2} (\tilde{\varphi}_y)^2 + \eta = 0,
\end{cases} 
\end{align*}   
where $B = -\alpha \gamma + \frac{1}{2} \alpha^2$ is the new Bernoulli's constant. Hence, $c = \gamma - \alpha$ is the true wave speed and 
$$
\gamma^2 + 2 B = \gamma^2 - 2 \alpha \gamma + \alpha^2 = (\gamma - \alpha)^2 = c^2,
$$
according to the definition of the bifurcation parameter $\gamma$ in Babenko's equation (\ref{BE}).
\end{proof}

By Proposition \ref{prop-comparison}, we need to normalize the speed appropriately due to the fact that the linear slope $\alpha$ is nonzero 
in the AFM method used in \cite{CDT22}. While doing so, we discovered 
a typo in \cite[Eq. (A.4e)]{CDT22} with the correct value:
$$
Q_{2,0} = \frac{1-c_0^4(h)}{4c_0^3(h)} = \frac{1}{2c_0(h) \sinh(2h)},
$$
such that $\alpha = a^2 Q_{2,0} + \mathcal{O}(a^4)$ according to 
\cite[Eq. (A.3b)]{CDT22}. By using the correction to the speed from \cite[Eq. (A.2b)]{CDT22} and the correction to the potential from 
the revised value for $Q_{2,0}$, we recover the correct expansion for the true wave speed:
\begin{align*}
c &= c_0(h_0) + \frac{a^2}{2 c_0(h_0)} \left[ 
\frac{9 - 6 \tanh^2(h_0) + 5 \tanh^4(h_0)}{8 \tanh^3(h_0)} - \frac{1}{2} \coth(h_0) {\rm sech}^2(h_0) \right]  + \mathcal{O}(a^4) \\
&= c_0(h_0) + \frac{a^2}{2 c_0(h_0)} \left[ 
\frac{9 - 10 \tanh^2(h_0) + 9 \tanh^4(h_0)}{8 \tanh^3(h_0)} \right]  + \mathcal{O}(a^4)
\end{align*}
in agreement with (\ref{exp-c}).  Thus, the Stokes expansions (\ref{eq:EtaExpansion}), (\ref{eq:c2Expansion}), (\ref{eq:BernouliExpansion}), and (\ref{mass}) obtained via the conformal transformation recovers the Stokes expansion from \cite{CDT22} obtained via the AFM formulation of Euler's equations.

\subsection{Comparison with \cite{Zakharov96}} 

The correction term (\ref{eq:c2Condition}) to the speed $c$ in terms of the conformal depth $h$ coincides with the expression in \cite[Eq. (4.13)]{Zakharov96}. The mean-zero constraint $M(\eta) = 0$ was enforced in \cite{Zakharov96} and the Bernoulli constant $B \neq 0$ was used in the relevant computations, see $b_0$ in \cite[Eqs. (4.5) and (4.9)]{Zakharov96}, although the details of computations were not explicitly written in \cite{Zakharov96}. Converting the expansion of the speed $c$ in terms of the physical depth $h_0$ yields (\ref{exp-c}), which coincides with the expression in \cite[Eq. (4.14)]{Zakharov96}. No integral constraint (\ref{eq:DepthIdentity-small}) was written explicitly in \cite{Zakharov96}.

\section{Numerical Experiments}\label{sec:numerics}

Here, we numerically compute solutions to Babenko's equation (\ref{eq:Babenko}) subject to the constraints (\ref{eq:DepthIdentity}) and (\ref{eq:Local}) for a 
variety of fixed physical depth $h_0 > 0$. The solutions are parameterized by the wave steepness $\mathcal{S}$, which is defined  as the difference between the maximal and minimal elevations in the profile $\eta \in C^{\infty}_{\rm per}$ divided by the wave period $2 \pi$, see \cite{DDS24,DLK16,DS23} for $h_0 = \infty$. The importance of this parameter is that the wave profile $\eta$ and other parameters of the system such as $\gamma$, $c$, $B$, and $h$ are 
smoothly parameterized in terms of $\mathcal{S}$. 

Due to the integral constraint (\ref{eq:DepthIdentity}), solutions 
of Babenko's equation (\ref{eq:Babenko}) differ from solutions of the classical Babenko equation (\ref{BE}), which have been studied numerically in a fixed conformal depth $h$ in \cite{SB25} and in a fixed physical depth $h_0$ in \cite{BDS26}. However, using (\ref{IC}) and the correspondence $c^2 = \gamma^2 + 2B(\gamma,h)$, we may compute the Bernoulli constant $B$ and the wave speed $c$ from the bifurcation parameter $\gamma$ used in \cite{BDS26,SB25} (there labelled as $c$). In Section \ref{sec:ComputingWaves} we discuss some improvements to the numerical methods of \cite{BDS26,SB25}, plot derivatives of Stokes waves with respect to $h$ and $c$, and provide comparisons with the asymptotics of Section \ref{sec-small}. In Section \ref{sec:constraints} 
we address the non-degeneracy and solvability conditions of Proposition 
\ref{prop-2-parameters} in the continuation of the solution family from the small-amplitude limit towards steep Stokes waves. 
In Section \ref{sec:LH+F} we recreate the results of Section 4 of \cite{BDS26} on the asymptotic behavior of the solution family for steep Stokes waves. 
Throughout this section, we set the mass constraint $M(\eta) = 0$. 

\subsection{Computing Stokes Waves} 
\label{sec:ComputingWaves}

\begin{figure}[htb!]
    \centering
    \includegraphics[width=\linewidth]{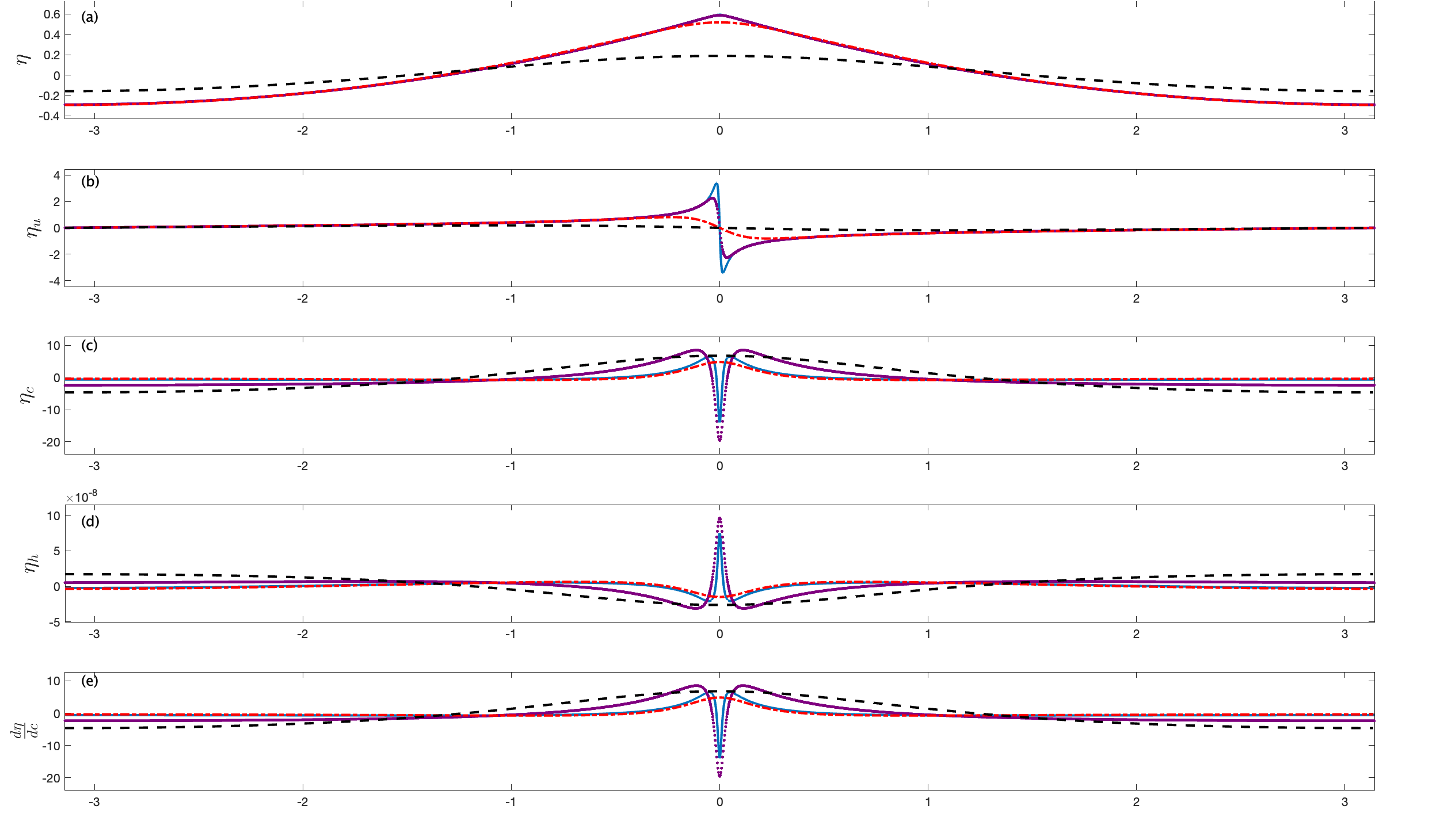}
    \caption{Profiles $\eta$ of Stokes waves of increasing steepness $\mathcal{S}$ versus the conformal variable $u$ (panel (a)), their partial derivatives in $u$ (panel (b)), $c$ (panel (c)), and $h$ (panel (d)), as well as their total derivative in speed $c$ (panel (e)), for $h_0 = 10$. In each plot, the solid blue line corresponds to the wave steepness $\mathcal{S}=0.14036$, the dotted purple line to $\mathcal{S}=0.13958$, the dash-dotted red line to $\mathcal{S}=0.12890,$ and the dashed black profile  to $\mathcal{S}=0.055011$.}
    \label{fig:EtaDerivs10}
\end{figure}

To solve the system (\ref{eq:Babenko}), (\ref{eq:DepthIdentity}), and (\ref{eq:Local}), we first solve Babenko's equation (\ref{BE}) using a Newton-Krylov method. Each step of the Newton--Krylov method is given by
\begin{align}
\label{eq:NewtonStep}
    \eta^{(N+1)}=\eta^{(N)}+\delta \eta^{(N)}, \quad 
    \delta \eta^{(N)}= -\mathcal{L}_h^{-1}S_h(\eta^{(N)}).
\end{align}
Once $\eta^{(N+1)}$ converges in the iterations of (\ref{eq:NewtonStep}), we have numerically computed a solution $\eta$ of Babenko's equation (\ref{BE}) with parameters $(\gamma,h)$. The Bernoulli constant $B = B(\gamma,h)$ is recovered from (\ref{IC}) and the speed $c$ is obtained from 
$c = \sqrt{\gamma^2 + 2B(\gamma,h)}$. 

Once the profile $\eta$ has been computed, we substitute $\eta$ into (\ref{eq:Local}) to compute the difference between $h - \overline{\eta}$ and the fixed level $h_0$ of the physical depth. This difference is used in a secant method to find a value of $h$ for which the equation (\ref{eq:Local}) is satisfied up to a maximum absolute error of $10^{-11}$.  We have then computed a solution of Babenko's equation (\ref{BE}) with parameters $(\gamma,h(\gamma))$ adjusted for the desired physical depth $h_0$. 

Note that the bifurcation parameter $\gamma$ in Babenko's equation (\ref{BE}) obtained in \cite{BDS26} is identical to the speed parameter in the AFM formulation obtained in \cite{CDT22}; see Proposition \ref{prop-comparison}. It is different from the correct wave speed $c$ in the original Babenko's equation (\ref{eq:Babenko}) for two reasons: it misses the Bernoulli constant in \cite{BDS26} and it includes the linear growth of the potential in \cite{CDT22}. 

\begin{figure}[htb!]
	\centering
	\includegraphics[width=\linewidth]{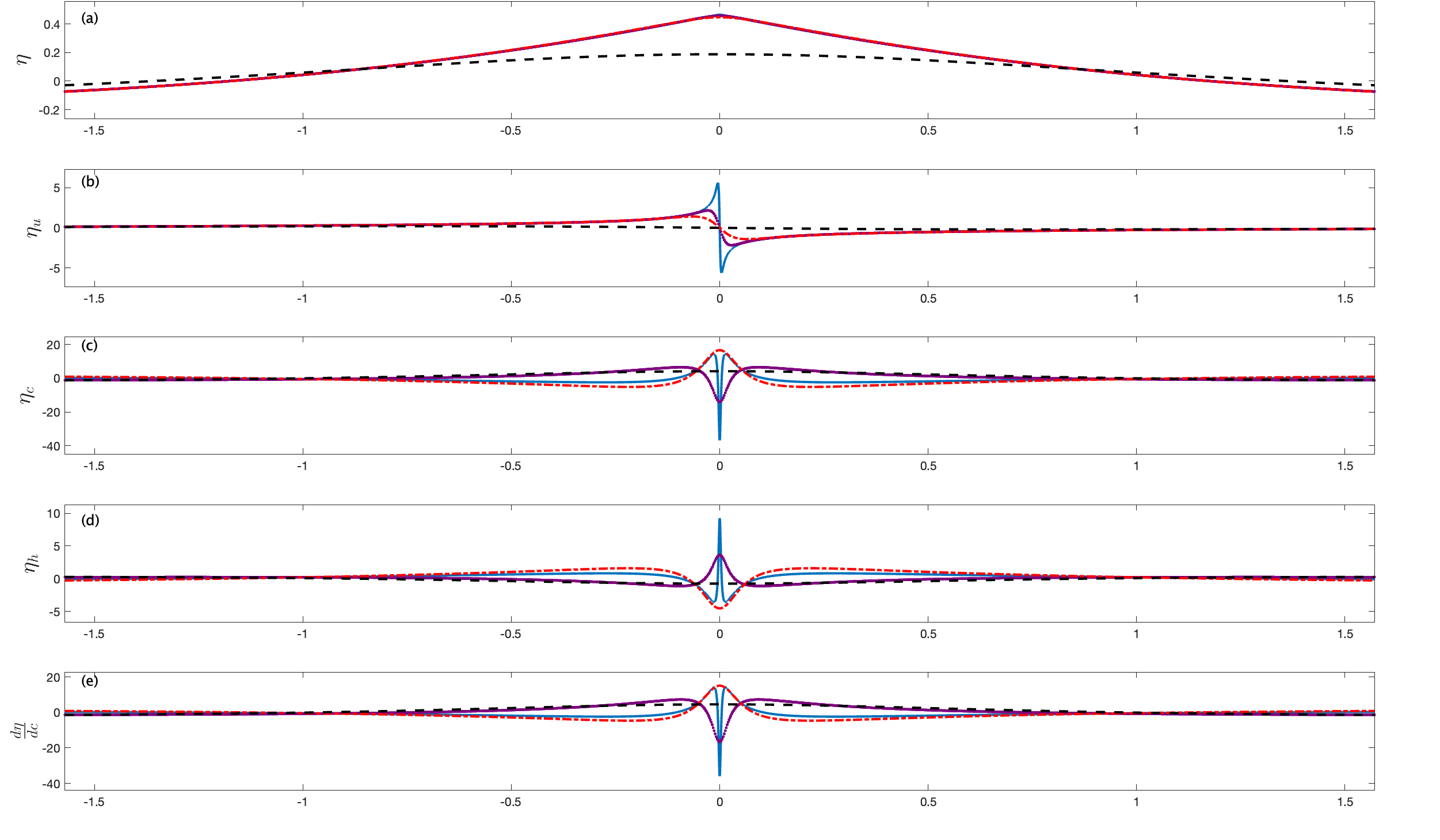}
	\caption{The same as Figure \ref{fig:EtaDerivs10} but for $h_0 = 1$. In each plot, the solid blue line corresponds to the wave steepness $\mathcal{S}=0.10034$, the dashed purple line to $\mathcal{S}=0.099289$, the dash-dotted red line to $\mathcal{S}=0.09770$, and the dotted black line to $\mathcal{S}=0.04934$. Plots are zoomed in on the middle half of a period of the wave to better resolve the structure near the crest.}
	\label{fig:EtaDerivs1}
\end{figure}

Having attained the profiles of Stokes waves with the wave speed $c$ at the physical depth $h_0$, we compute the partial derivatives of their profiles with respect to $u$, $c$, and $h$ under the assumptions of Proposition \ref{prop-2-parameters}: ${\rm Ker}(\mathcal{L}_h) = {\rm span}(\eta')$ and $1 - \overline{\partial_h \eta} \neq 0$. Taking a derivative with respect to $u$ is direct as it is diagonal in Fourier space. Taking the derivatives with respect to $h$ and $c$ is more complicated, as they require the solution of the systems (\ref{derivative-2})--(\ref{derivative-3}) and  (\ref{derivative-Ber-1})--(\ref{derivative-Ber-2}). We obtain $\partial_hB$ and $\partial_cB$ from these systems in the quotient form:
\begin{align}
\label{eq:BhSoln}
    \partial_hB =\frac{2\langle\mathcal{L}_h^{-1}\partial_hS_h(\eta), \partial_hS_h(\eta)\rangle-\langle (c^2-2B)\eta-\eta^2,(\partial_h^2K_h)\eta\rangle}{4\langle \mathcal{L}_h^{-1}K_h\eta,\partial_hS_h(\eta)\rangle-2\langle \eta,(\partial_hK_h)\eta\rangle-4\pi}
\end{align}
and
\begin{align}
    \label{eq:BcSoln}
    \partial_cB = c\frac{2\langle \mathcal{L}_h^{-1}K_h\eta,\partial_hS_h(\eta)\rangle-\langle \eta,(\partial_hK_h)\eta\rangle}{2\langle \mathcal{L}_h^{-1}K_h\eta,\partial_hS_h(\eta)\rangle-\langle \eta,(\partial_hK_h)\eta\rangle-2\pi}.
\end{align}
Substituting (\ref{eq:BhSoln}) and (\ref{eq:BcSoln}) into the linear inhomogeneous 
equations (\ref{derivative-2}) and (\ref{derivative-3}) and numerically inverting the linear operator $\mathcal{L}_h$ yields expressions for $\partial_h\eta$ and $\partial_c\eta$. The total derivative of the surface elevation $\eta$ with respect to the wave speed $c$ is computed from (\ref{derivative-constraints}) by the chain rule:
\begin{align*}
\frac{d\eta}{dc} = \partial_c\eta+\frac{dh}{dc}\partial_h\eta = \partial_c\eta+\frac{\overline{\partial_c\eta}}{1-\overline{\partial_h\eta}} \partial_h\eta.
\end{align*}

\begin{figure}[htb!]
    \centering
    \includegraphics[width=\linewidth]{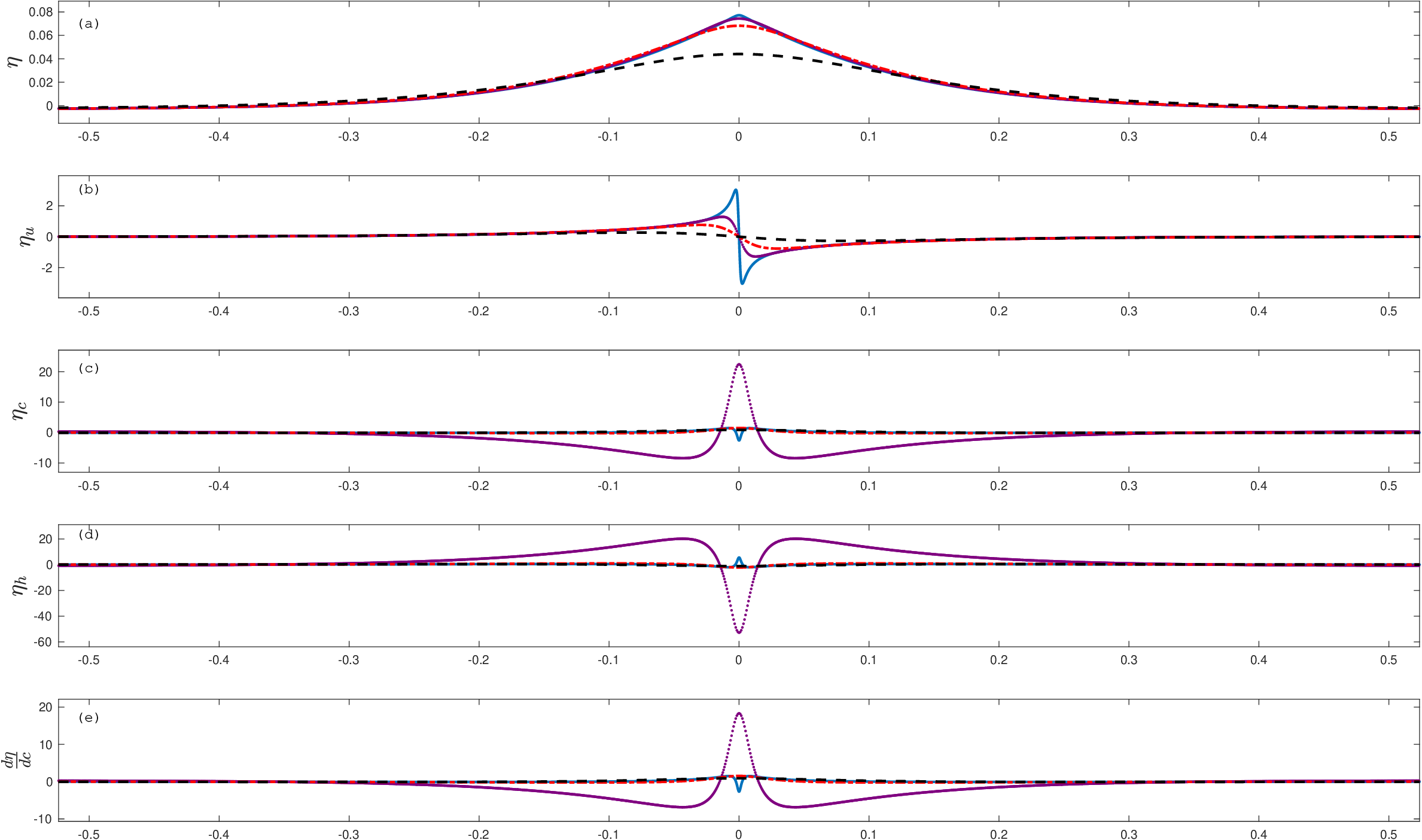}
    \caption{The same as Figure \ref{fig:EtaDerivs10} but for $h_0 = 0.1$. In each plot, the solid blue line corresponds to the wave steepness $\mathcal{S}=0.012757$, the dashed purple line to $\mathcal{S}=0.012312$, the dash-dotted red line to $\mathcal{S}=0.0113446$ , and the dotted black line to $\mathcal{S}=0.007441$. Plots are zoomed in on the middle sixth of a period of the wave to better resolve the structure near the crest.}
    \label{fig:EtaDerivs01}
\end{figure}

Figures \ref{fig:EtaDerivs10},  \ref{fig:EtaDerivs1}, and \ref{fig:EtaDerivs01} present profiles of $\eta$, $\eta'$, $\partial_c \eta$, $\partial_h \eta$, and $\frac{d\eta}{dc}$ for $h_0 = 10$, $h_0 = 1$, and $h_0 = 0.1$, respectively. In deep water, the profile of $\frac{d\eta}{dc}$ largely resembles the profile of $\partial_c \eta$, but as depth $h_0$ decreases, they begin to incorporate more of the profile of $\partial_h \eta$. 

In the small-amplitude regime, the parameter $\gamma$, the wave speed $c$, and the wave steepness $\mathcal{S}$ are increased together. However, for the waves of increasing steepness $\mathcal{S}$, there exist bifurcation points, where 
both the bifurcation parameter $\gamma$ and the wave speed $c$ reach a first maximum. Past the bifurcation points (which are different for $\gamma$ and $c$), $\gamma$ and $c$ decrease with increasing $\mathcal{S}$ until they reach the second such bifurcation points, where $\gamma$ and $c$ are minimal. These oscillations are discussed further in Section \ref{sec:LH+F}.  

\begin{figure}[htb!]
    \centering
    \includegraphics[width=\linewidth]{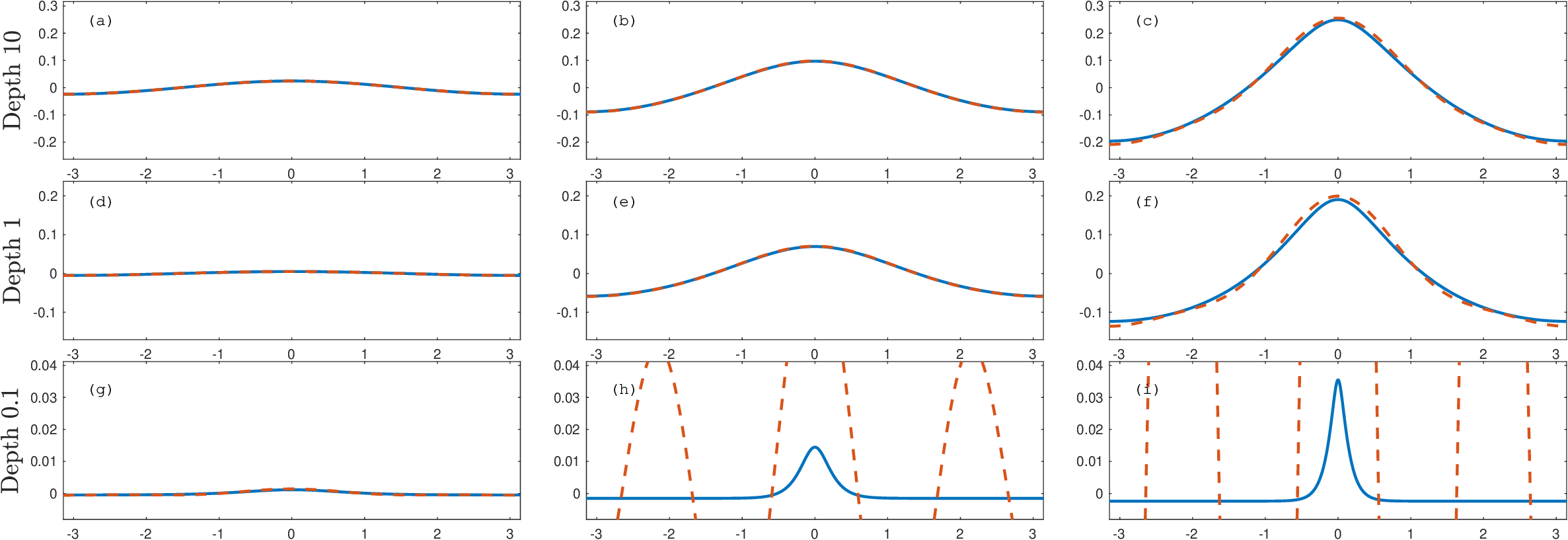}
    \caption{Comparison between numerically computed profiles $\eta$ (solid blue line) and those computed from the expansion (\ref{eq:EtaExpansion}) (dashed red line) in depths $h_0 = 10$, 
    	$h_0 = 1$, and $h_0 = 0.1$ (panels (a-c), (d-f), and (g-i)) respectively. Panels (a), (d), and (g) each show a wave of about 2 percent of the amplitude of the steepest wave, with amplitude parameters $a = 0.0242,0.005019,$ and $7.688\cdot 10^{-4}$ respectively. Panels (b), (e), and (h) each show a wave of about 20 percent of the amplitude of the steepest wave, with amplitude parameters $a = 0.09242, 0.06346,$ and $0.00506$ respectively. Panels (c), (f), and (i) each show a wave of about 50 percent of the amplitude of the steepest wave and have amplitude parameters $a = 0.2167, 0.15015,$ and $0.00948$ respectively.}
    \label{fig:stokesCompare}
\end{figure}

\begin{figure}[htb!]
	\centering
	\includegraphics[width=\linewidth,height=0.47\textwidth]{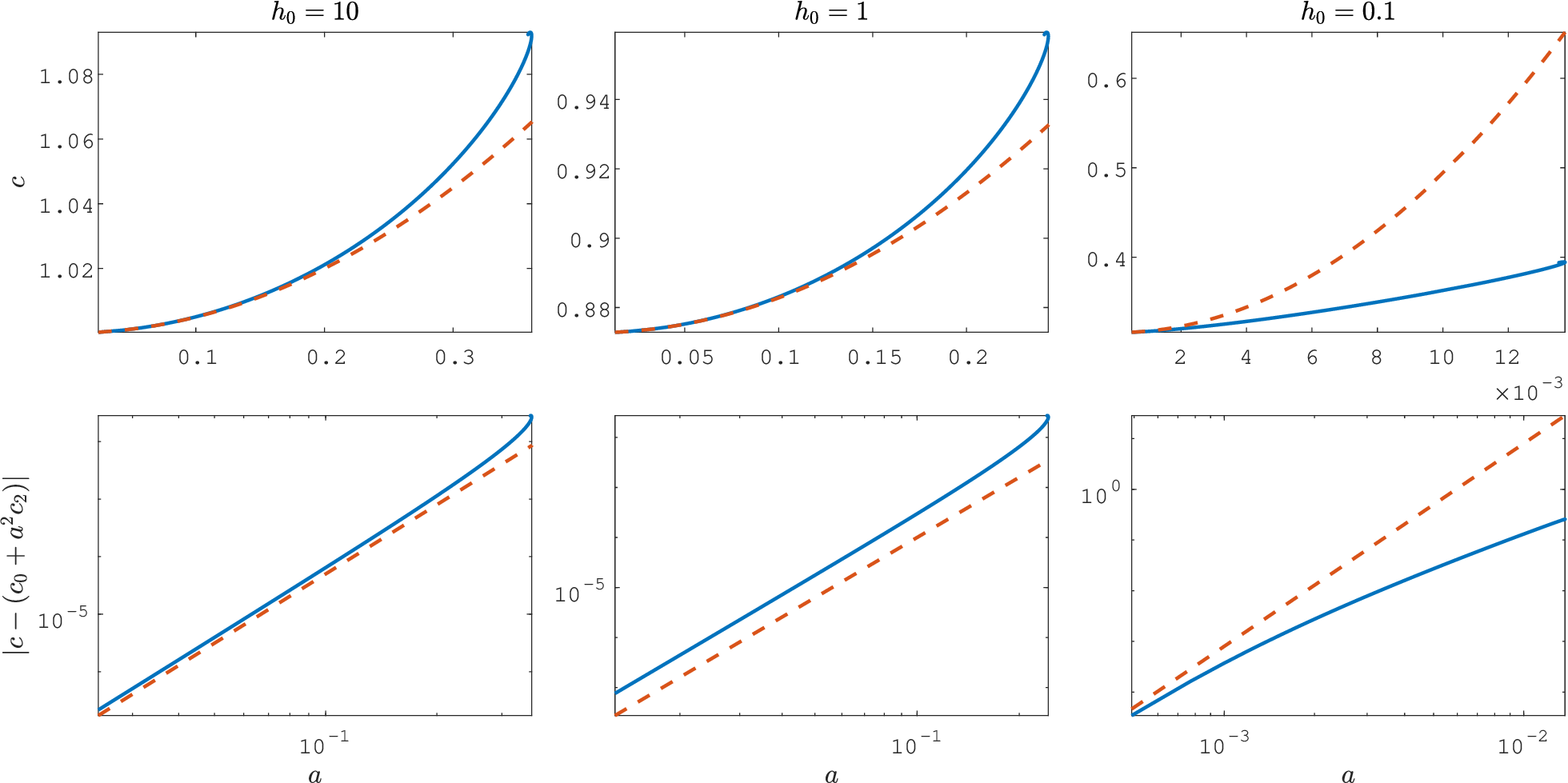}
	\caption{Top: Comparison between the asymptotic approximation of $c_0 + a^2 c_2$ given by  (\ref{exp-c}) (dashed red lines) and the numerically computed value for the wave speed $c$ (solid blue lines) in depths $h_0 = 10$ (left), $h_0 = 1$ (middle), and $h_0 = 0.1$ (right). Bottom: residual $c-(c_0+a^2c_2)$ (solid blue line) and plot of $Ca^4$ to demonstrate order of error term (dashed red line), on log--log axes.}
	\label{fig:speedCompare}
\end{figure}

\begin{figure}[htb!]
	\centering
	\includegraphics[width=\linewidth,height=0.47\textwidth]{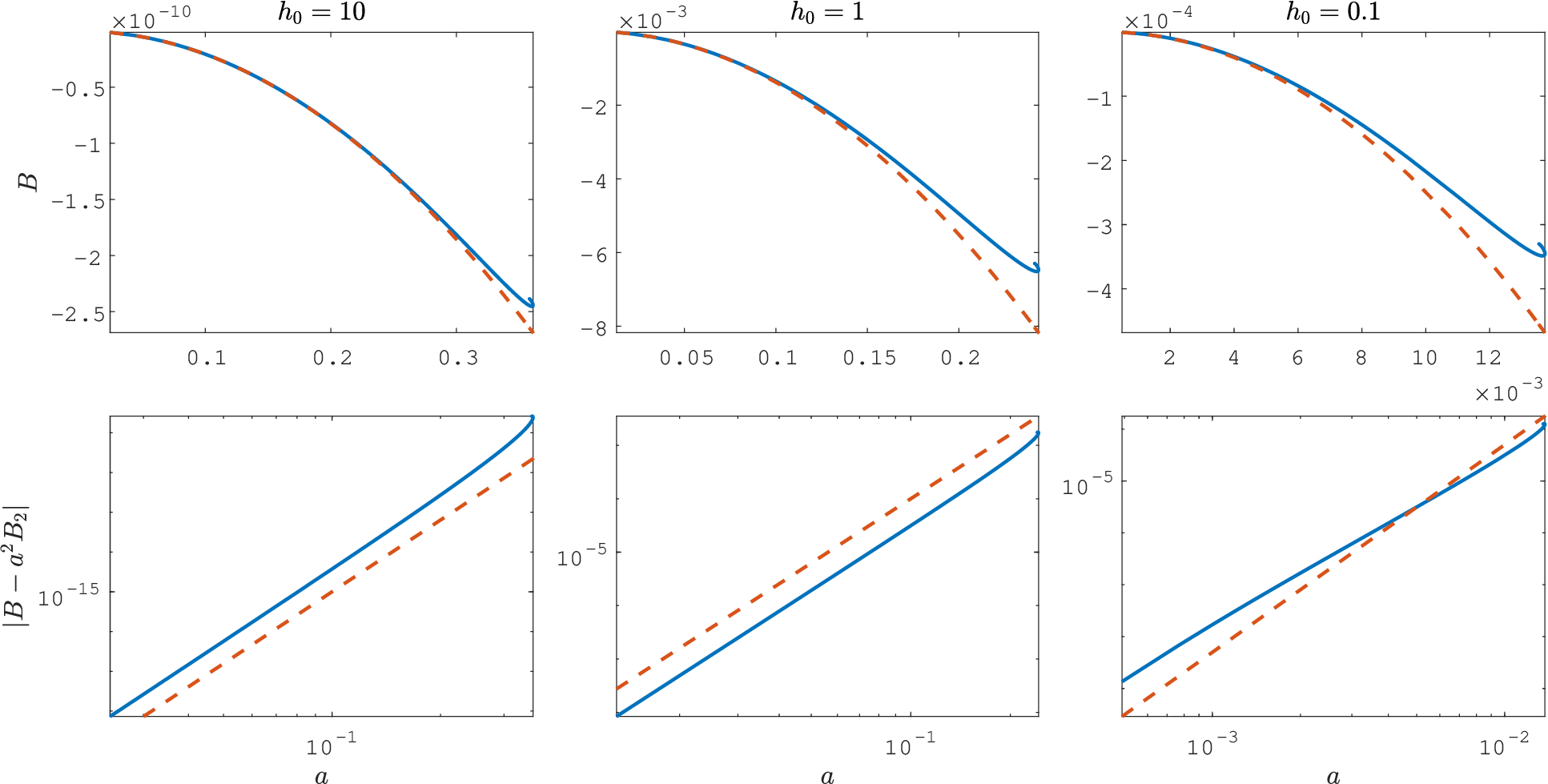}
	\caption{Top: Comparison between the asymptotic approximation of $a^2 B_2$ given by (\ref{exp-B}) (dashed red lines) and the numerically computed value for Bernoulli's constant $B$ (solid blue lines) in depths $h_0 = 10$ (left), $h_0 = 1$ (middle), and $h_0 = 0.1$ (right). Bottom: residual $B-a^2B_2$ (solid blue line) and plot of $Ca^4$ to demonstrate order of error term (dashed red line), on log--log axes.}
	\label{fig:bernoulliCompare}
\end{figure}

Figures \ref{fig:stokesCompare}, \ref{fig:speedCompare}, and \ref{fig:bernoulliCompare} compare respectively the surface elevation $\eta$, the wave speed $c$, and the Bernoulli constant $B$ with the Stokes expansions (\ref{eq:EtaExpansion}), (\ref{eq:c2Expansion}), and (\ref{eq:BernouliExpansion}). Instead of the wave steepness $\mathcal{S}$, we can use the wave amplitude $a$ defined in Section \ref{sec-small} to parameterize the solution family in the small-amplitude regime.

There is excellent agreement between the numerics and the asymptotics for small $a$. However, as the amplitude $a$ increases, and especially as the depth $h_0$ decreases, the quality of the match decreases. In the depth $h_0 = 0.1$ specifically, the asymptotic series appears to work for much smaller values of $a$, with next to no agreement in the data in panels (h) and (i). This is related to the growth of the coefficients of the expansion (\ref{eq:EtaExpansion}) as depth $h_0$ decreases, with values as large as $A_2=7.70367\cdot10^2$ and $A_3=4.5044\cdot10^5$ in depth $h_0 = 0.1$. For well ordering of the expansion in this depth $h_0$, $a$ needs to be less than $10^{-3}$.

\subsection{Continuation of the solution family}
\label{sec:constraints}

Computations of the solution family by using the Newton--Krylov method in Section \ref{sec:ComputingWaves} and  the partial derivatives of the wave profile $\eta$ with respect to $c$ and $h$ are only possible under the assumptions of Proposition \ref{prop-2-parameters}: ${\rm Ker}(\mathcal{L}_h) = {\rm span}(\eta')$ and $1 - \overline{\partial_h \eta} \neq 0$.
We have found numerically that the non-degeneracy assumption ${\rm Ker}(\mathcal{L}_h) = {\rm span}(\eta')$ is not satisfied at each critical point of the bifurcation parameter $\gamma$ for a fixed \textit{conformal} depth $h$ and $1 - \overline{\partial_h \eta} \neq 0$ is not satisfied at each critical point of the bifurcation parameter $\gamma$ for a fixed \textit{physical} depth $h_0$, when $\gamma$ is plotted versus steepness $\mathcal{S}$.

\begin{figure}[htb!]
	\centering
	\includegraphics[width=\linewidth,height=0.47\textwidth]{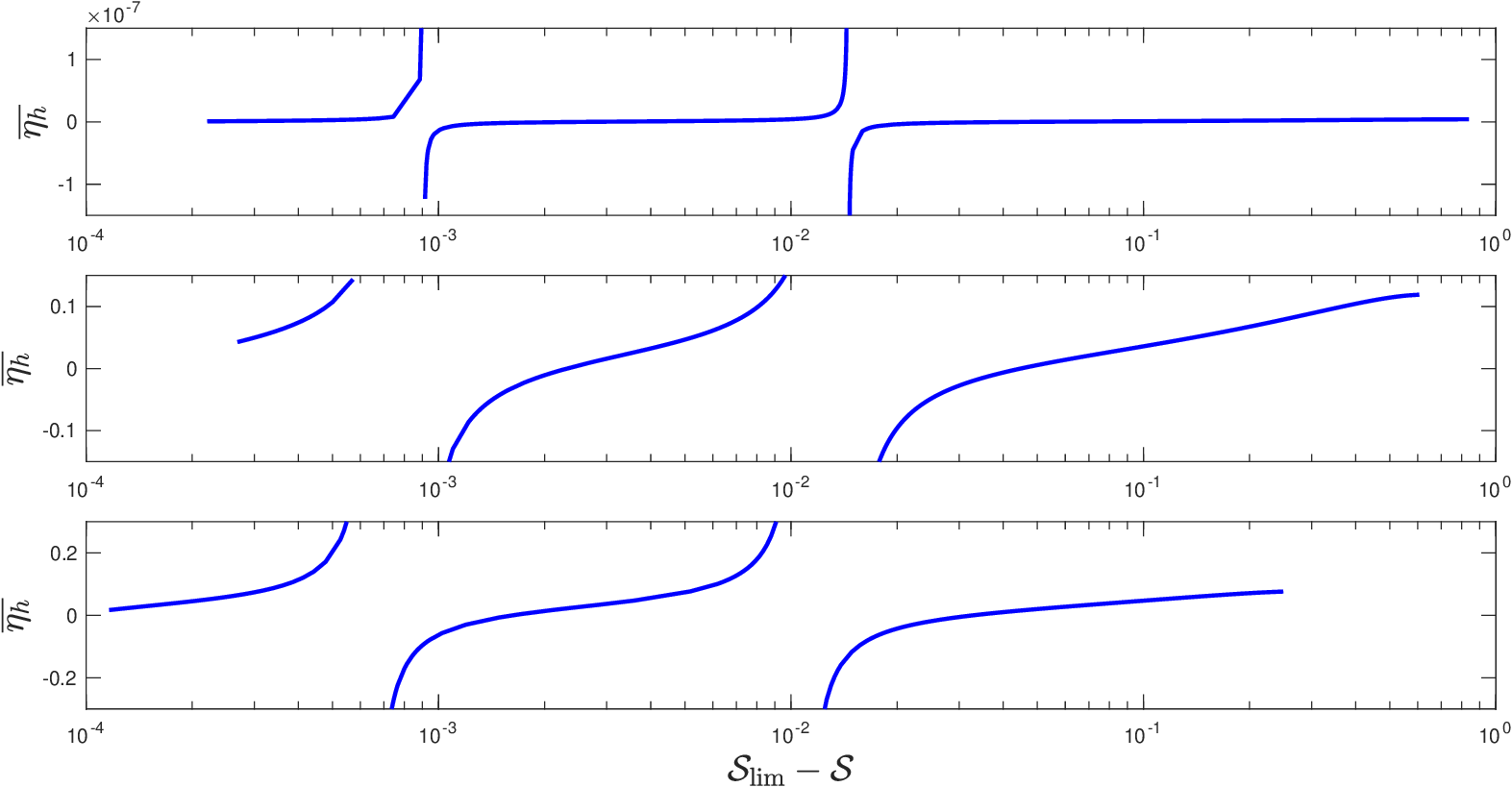}
	\caption{Plots of the numerically computed  $\overline{\partial_h \eta}$ for $h_0 = 10$ (top), $h_0 = 1$ (middle), and $h_0 = 0.5$ (bottom). Each is plotted on a semilog scale with the logarithmic horizontal axis showing the difference in steepness $\mathcal{S}$ from the limiting wave steepness $\mathcal{S}_{\rm lim}$. The divergence of $\overline{\partial_h \eta}$ occurs when the dimension of the Jacobian $\mathcal{L}_h$ increases to two}
	\label{fig:avgetahFits}
\end{figure}

Figure \ref{fig:avgetahFits} displays the dependence of $\overline{\partial_h \eta}$ on the relative steepness $\mathcal{S}_{\rm lim} - \mathcal{S}$ from the limiting wave. The small-amplitude regime corresponds to the right side of this figure and the regime of steeper Stokes waves corresponds to the left side. We observe that $\overline{\partial_h \eta}$ diverges to $-\infty$ before the dimension of the Jacobian $\mathcal{L}_h$ increases and diverges to $+\infty$ after this point. Therefore, the solvability condition $1 - \overline{\partial_h \eta} \neq 0$ fails past each such point as $\mathcal{S}$ increases. Assuming a smooth parameterization in steepness $\mathcal{S}$, we have for a fixed physical depth $h_0$ that 
\begin{align} 
\label{eq:constraintDerivative}
    0 = (1-\overline{\partial_h \eta}) \frac{dh}{d\mathcal{S}} -\overline{\partial_\gamma \eta} \frac{d\gamma}{d\mathcal{S}}.
\end{align}
If $\overline{\partial_h \eta}=1$, $\overline{\partial_{\gamma} \eta} \neq 0$, and $\frac{dh}{d\mathcal{S}}$ is finite, 
then it follows from (\ref{eq:constraintDerivative}) that $\frac{d\gamma}{d\mathcal{S}}=0$ so that the point 
where $\overline{\partial_h \eta} = 1$ corresponds to the critical point of the bifurcation parameter $\gamma$ 
plotted versus steepness $\mathcal{S}$ for a fixed physical depth $h_0$. Numerical evidence 
supports this conjecture and show that $h$ varies smoothly across the points where $\overline{\partial_h \eta} = 1$ 
as a function of $\mathcal{S}$ and that $\overline{\partial_{\gamma} \eta}$ is bounded away from zero near these points. 

On the other hand, $\overline{\partial_h \eta}$ diverges at each point where the dimension of ${\rm Ker}(\mathcal{L}_h)$ increases to $2$. 
We have verified numerically that these points correspond to the critical points of $\gamma$ with respect to steepness 
$\mathcal{S}$ for a fixed conformal depth $h$. That is, ${\rm Ker}(\mathcal{L}_h)$ is two-dimensional precisely when 
$\frac{\partial \gamma}{\partial \mathcal{S}}=0$ compared to 
$$
\frac{d \gamma}{d\mathcal{S}} = \frac{\partial \gamma}{\partial \mathcal{S}} + \frac{\partial \gamma}{\partial h} \frac{dh}{\mathcal{S}},
$$
which is valid for a fixed physical depth $h_0$. This criterion corresponds to the numerical results of \cite{SB25}, 
which have suggested that the Jacobian of Babenko's equation (\ref{BE}) has a double-zero eigenvalue at critical points of 
$\gamma$ for a fixed conformal depth $h$.

We note that the divergence of $\overline{\partial_h \eta}$ has a strong impact on the convergence of the secant method in finding solutions of $h = h_0 + \overline{\eta}$, as the derivative of $h_0$ with respect to $h$ diverges, decreasing the stability of the secant method and forcing us to take smaller step sizes, so that approaching the critical points of $\gamma$ versus $\mathcal{S}$ from below becomes prohibitive.

\begin{figure}[htb!]
	\centering
	\includegraphics[width=\linewidth,height=0.47\textwidth]{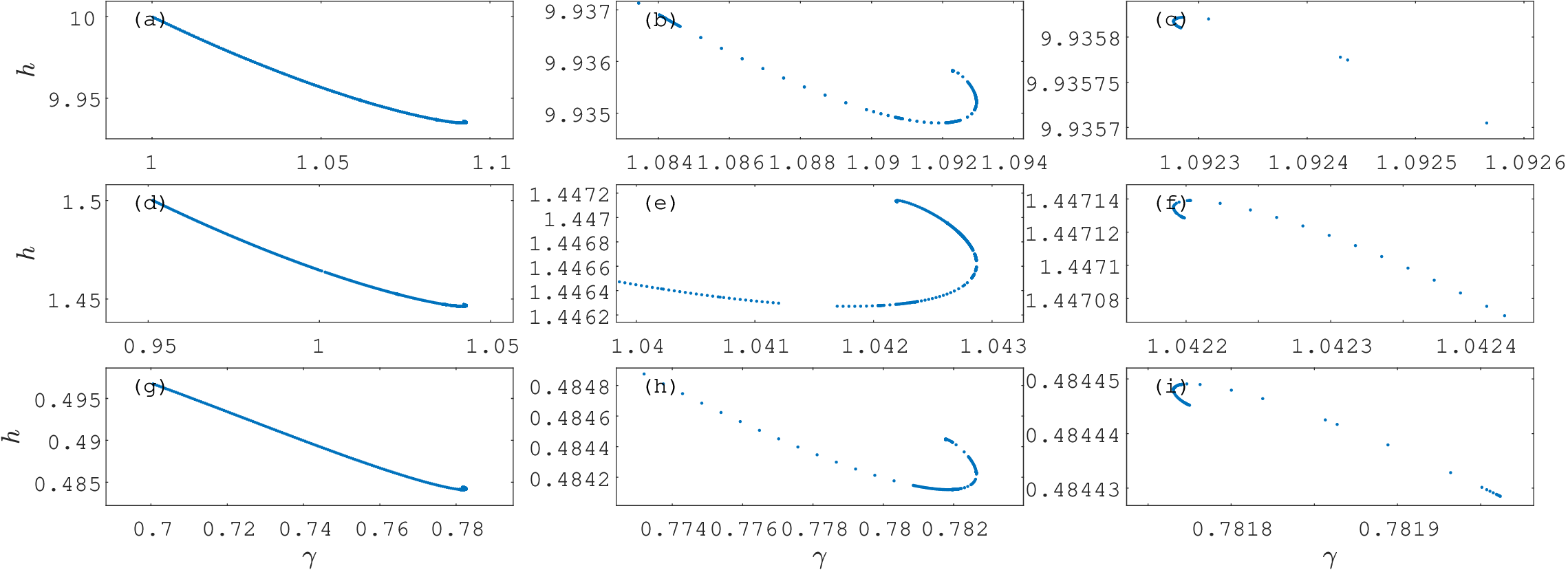}
	\caption{Plots of the conformal depth $h$ as a function of the bifurcation parameter $\gamma$ for $h_0 = 10$ (a-c), $h_0 = 1.0$ (d-f), and $h_0 = 0.5$ (g-i). The middle column is a zoom of the left column and the right column is a zoom of the middle column.}
	\label{fig:speedh}
\end{figure}

We also plot the conformal depth $h$ versus the bifurcation parameter $\gamma$ in Figure \ref{fig:speedh} to illustrate the spiral behavior of both $h$ and $\gamma$ with a smooth parameterization by steepness $\mathcal{S}$.

\subsection{Oscillations for steep Stokes waves}
\label{sec:LH+F}

Here, we report on the oscillations of the wave speed $c$, the Bernoulli constant $B$, and the wave energy $\mathcal{H} = \mathcal{H}(c,h(c))$, see (\ref{energy-wave}), as the steepness $\mathcal{S}$ is increased towards the limiting Stokes wave with the peaked profile. It is well-known from the asymptotic theory \cite{LHF1,LHF2} and numerical approximations \cite{DH26} for $h_0 = \infty$ that there are infinitely many oscillations, but due to the exponentially small phenomena near the limiting wave, we are only illustrating the first two oscillations.

The oscillation phenomena were observed for solutions of the classical Babenko equation (\ref{BE}) for $h_0 < \infty$ in \cite{BDS26}. As discussed in Section \ref{sec:ComputingWaves}, wave profiles found by solving Babenko's equation (\ref{eq:Babenko}) with constraints (\ref{eq:DepthIdentity}) and (\ref{eq:Local}) are identical to those from solving (\ref{BE}), with a slightly different value of the speed $c = \sqrt{\gamma^2 + 2B(\gamma,h)}$ due to the Bernoulli constant $B =  B(\gamma,h)$. 

\begin{figure}[htb!]
	\centering
	\includegraphics[width=\linewidth,height=0.47\textwidth]{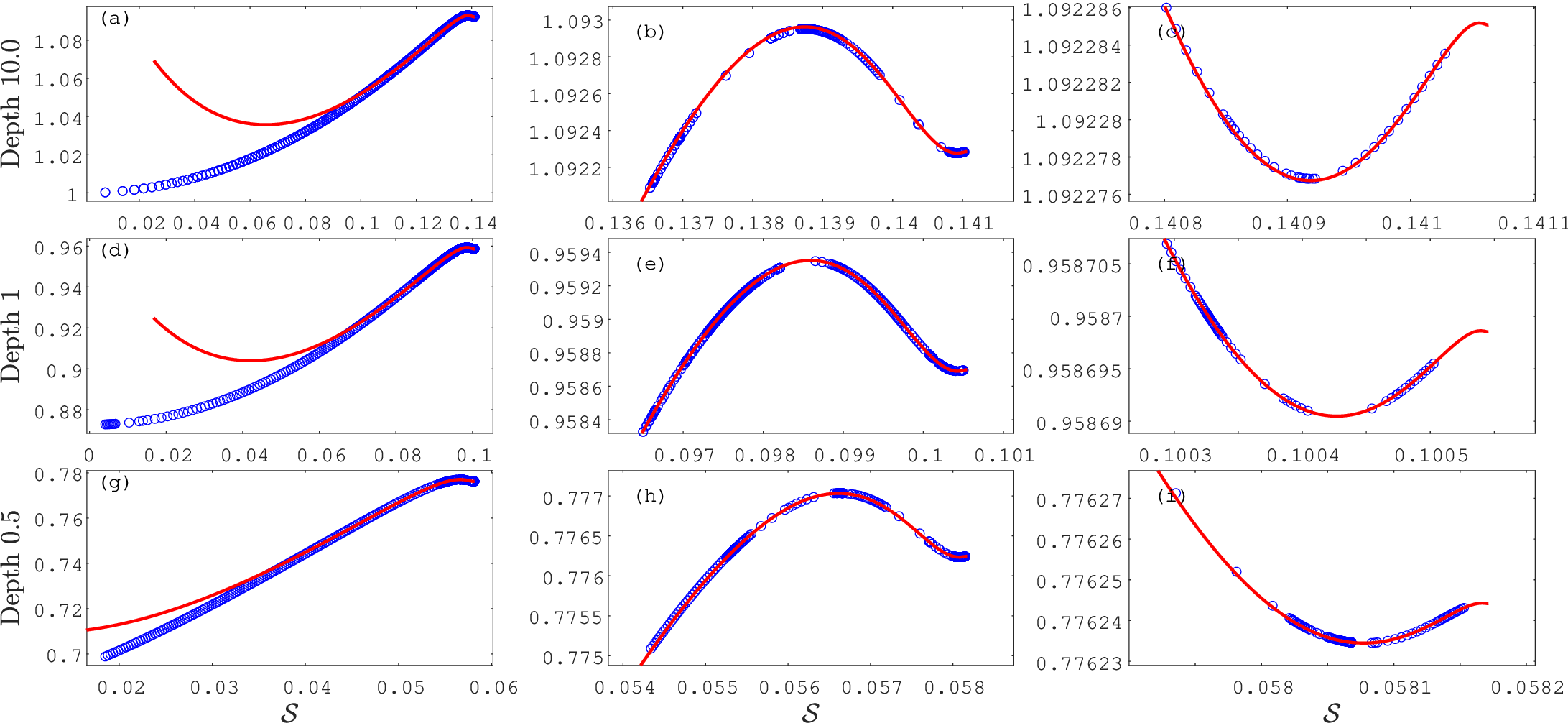}
	\caption{Wave speed $c$ of the numerically computed Stokes waves versus steepness $\mathcal{S}$ for $h_0 = 10$ (a-c), $h_0 = 1$ (d-f), and 
		$h_0 = 0.5$ (g-i). Circular data points are the numerically computed data and the solid line plots the fit from the asymptotic formula (\ref{eq:LHFoxAsymptotics}). The middle column is a zoom of the left column and the right column is a zoom of the middle column.}
	\label{fig:speedFits}
\end{figure}

\begin{figure}[htb!]
	\centering
	\includegraphics[width=\linewidth,height=0.47\textwidth]{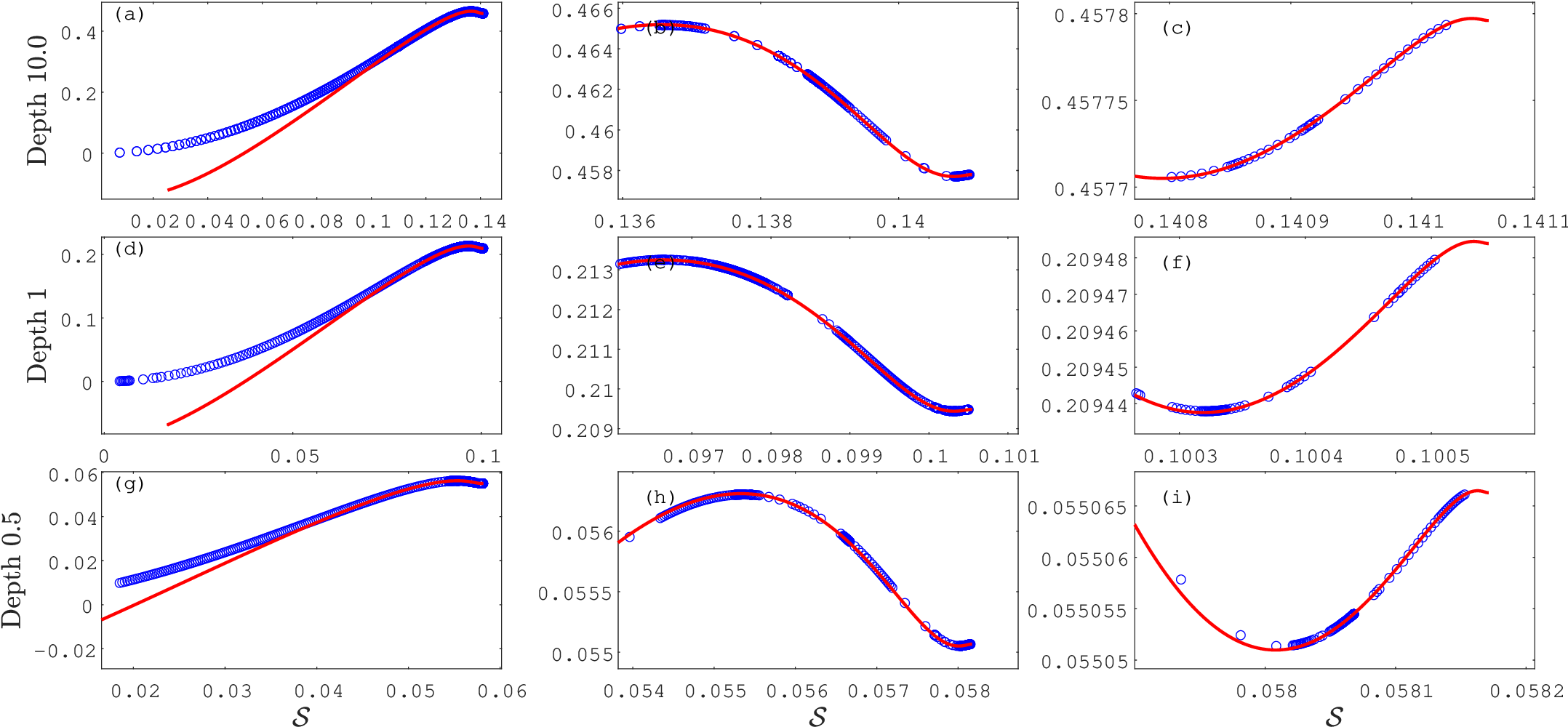}
	\caption{The same as Figure \ref{fig:speedFits} but fitting the  wave energy $\mathcal{H}$ to LHF's asymptotics and plotting versus steepness $\mathcal{S}$.}
	\label{fig:energyFits}
\end{figure}

\begin{figure}[htb!]
    \centering
    \includegraphics[width=\linewidth,height=0.47\textwidth]{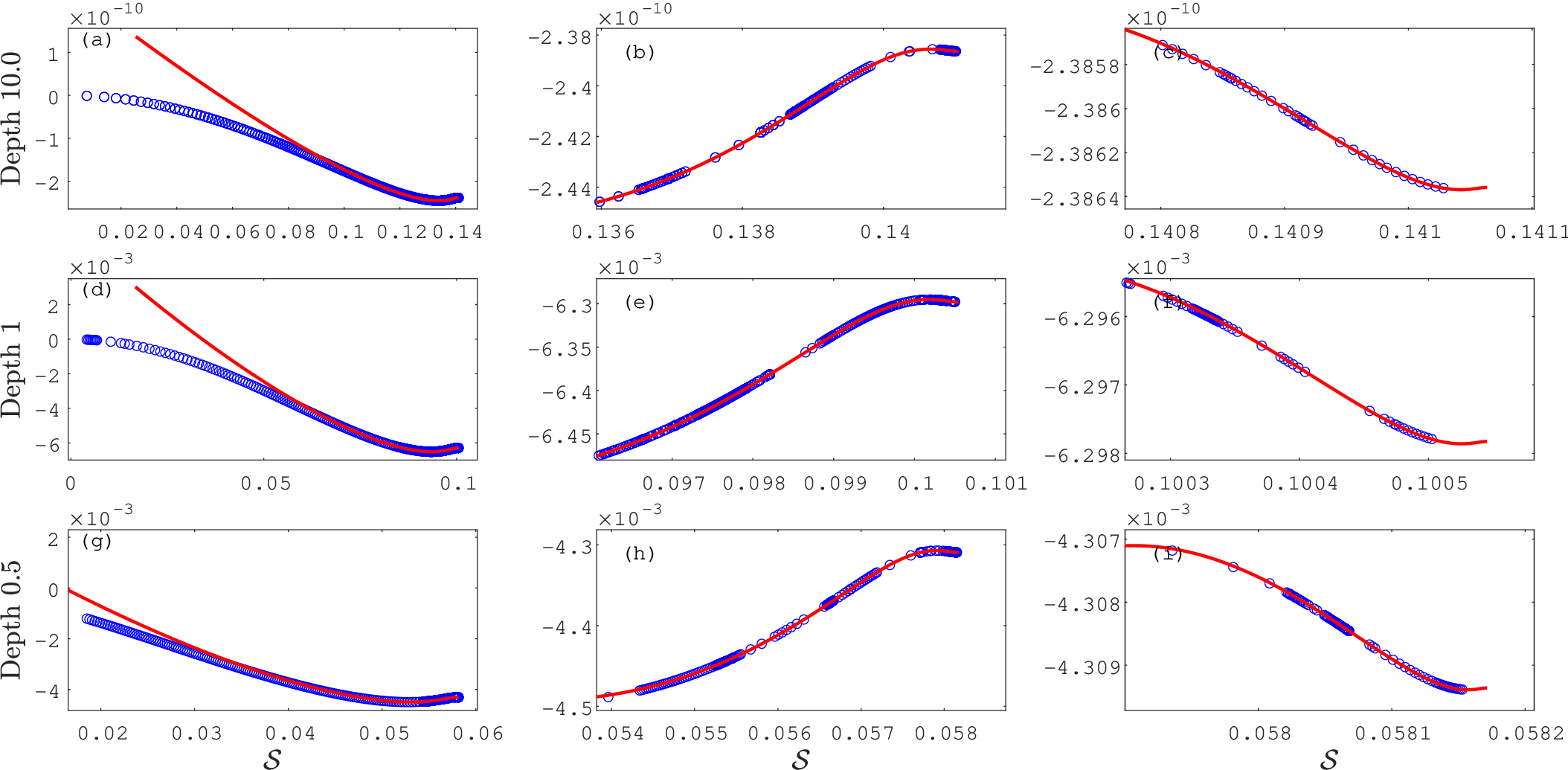}
    \caption{The same as Figure \ref{fig:speedFits} but fitting the  Bernoulli constant $B$ to LHF's asymptotics and plotting versus steepness $\mathcal{S}$.}
    \label{fig:bernoulliFits}
\end{figure}

To investigate if the oscillation phenomena holds for the system  (\ref{eq:Babenko}), (\ref{eq:DepthIdentity}), and (\ref{eq:Local}), we plot the changes in the speed $c$, the energy $\mathcal{H}$, and the Bernoulli constant $B$ versus steepness $\mathcal{S}$ in Figures \ref{fig:speedFits}, \ref{fig:energyFits}, and \ref{fig:bernoulliFits}, respectively. These figures demonstrate that the speed $c$ and energy $\mathcal{H}$ still oscillate in $\mathcal{S}$ for a variety of depths $h_0$. 

We superpose the numerical data in Figures \ref{fig:speedFits}, \ref{fig:energyFits}, and \ref{fig:bernoulliFits} with the outcomes of the asymptotic theory developed by Longuet-Higgins and Fox (LHF) in the 1980s \cite{LHF1,LHF2}. The correspondence is remarkably close. The LHF asymptotic results were written in terms of a perturbation parameter 
\begin{align*}
    \epsilon=\frac{s}{c_0\sqrt{2}},
\end{align*}
where $s = \sqrt{\gamma^2-2\eta(0)}$ is the speed of the Stokes wave at the crest and $c_0 = \sqrt{\tanh(h_0)}$. 
The perturbation parameter $\epsilon$ is proportional to the radius of curvature of the Stokes wave at the crest. As the limiting steepest wave with the peaked profile is approached, the radius of curvature and thus $\epsilon$ tend to zero. 
To create the plots shown in Figures \ref{fig:speedFits}, \ref{fig:energyFits}, and \ref{fig:bernoulliFits} we fit our numerical data to the asymptotic formula
\begin{align}
\label{eq:LHFoxAsymptotics}
    c^2\sim c_{\lim}^2+a_2\epsilon^3\cos(\mu\ln(\epsilon)+a_3),
\end{align}
where $\mu \approx 2.143$ is a parameter from the theory \cite{LHF1,LHF2}. $c_{\lim}, a_2$, and $a_3$ are then chosen to minimize the least squares error between the asymptotic formula \eqref{eq:LHFoxAsymptotics}, weighting the data with a factor of $1/\epsilon^6$ to account for the $\mathcal{O}(\epsilon^6)$ next order term. 

\begin{figure}[htb!]
	\centering
	\includegraphics[width=\linewidth,height=0.47\textwidth]{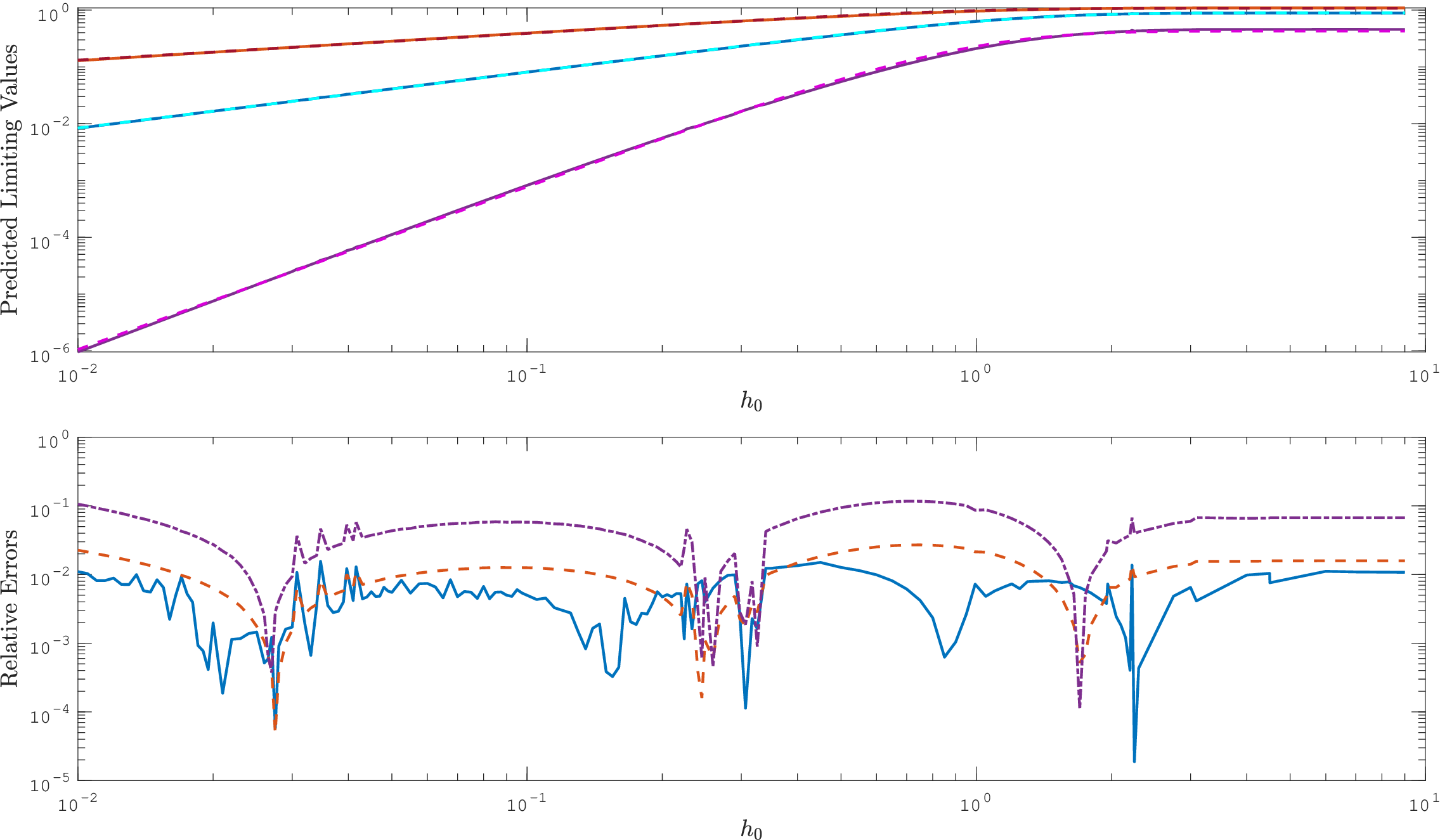}
	\caption{Top: the predicted limiting wave speed (red) and energy (purple) of the steepest wave as a function of depth, as well as fits to powers of the dispersion relation. Also shown in blue is the steepness of the limiting Stokes wave. Bottom: the relative residuals of these fits, which demonstrate that the errors are consistently an order of magnitude smaller than the data.}
	\label{fig:extremeSpeedsEnergies}
\end{figure}

These asymptotics tell us two things. First, since the $\epsilon\to0$ limit describes the steepest wave, $c_{\lim}^2$ is appropriately the squared speed of the steepest wave. Second, they tell us that the largest oscillation of a quantity described by these asymptotics is the first one. Thus, the first maximum of that quantity is the largest in magnitude. Figure \ref{fig:extremeSpeedsEnergies} shows the limiting wave speed $c$ and energy $\mathcal{H}$ computed from (\ref{eq:LHFoxAsymptotics}) numerically.
We observe that the two quantities decay like powers of the dispersion relation $\omega^2(kh)=kh\tanh(kh)$ as in \cite{BDS26}. 

\begin{figure}[htb!]
    \centering
    \includegraphics[width=\linewidth,height=0.47\textwidth]{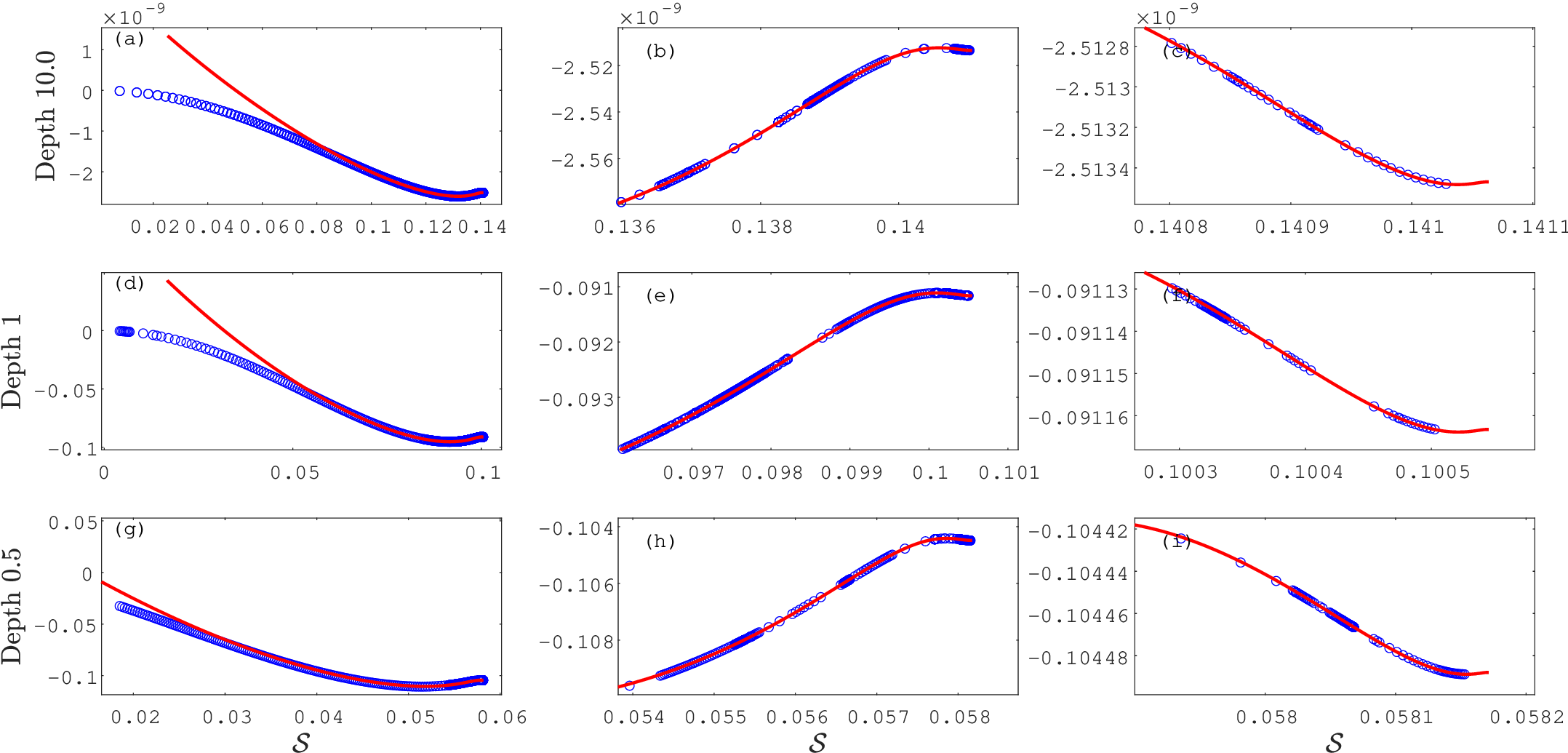}
    \caption{The numerically computed constant $\langle \eta,(\partial_h K_h)\eta\rangle$ versus steepness $\mathcal{S}$ for $h_0 = 10$ (a-c), $h_0 = 1.0$ (d-f), and $h_0 = 0.5$ (g-i). Circular data points are the numerically computed data and the solid line plots the fit from the asymptotic formula (\ref{eq:LHFoxAsymptotics}). The middle column is a zoom of the left column and the right column is a zoom of the middle column.}
    \label{fig:etaKetaFits}
\end{figure}

Finally, we address the constraint $\langle \eta,(\partial_hK_h)\eta\rangle \neq - 2\pi$ used in Theorem \ref{theorem-multiplicity}.
Figure \ref{fig:etaKetaFits} gives the plot of $\langle \eta,(\partial_hK_h)\eta\rangle$ versus steepness $\mathcal{S}$. 
The numerical constant $\langle \eta,(\partial_hK_h)\eta\rangle$ oscillates according to LHF's asymptotic formula (\ref{eq:LHFoxAsymptotics}).  The match between this formula and our data allows us to use them to bound the maximum value of $\langle \eta,(\partial_h K_h)\eta\rangle$ far away from $-2\pi$ providing evidence for the constraint (\ref{IFT-condition}) to hold everywhere. 

\begin{figure}[htb!]
    \centering
    \includegraphics[width=\linewidth,height=0.47\textwidth]{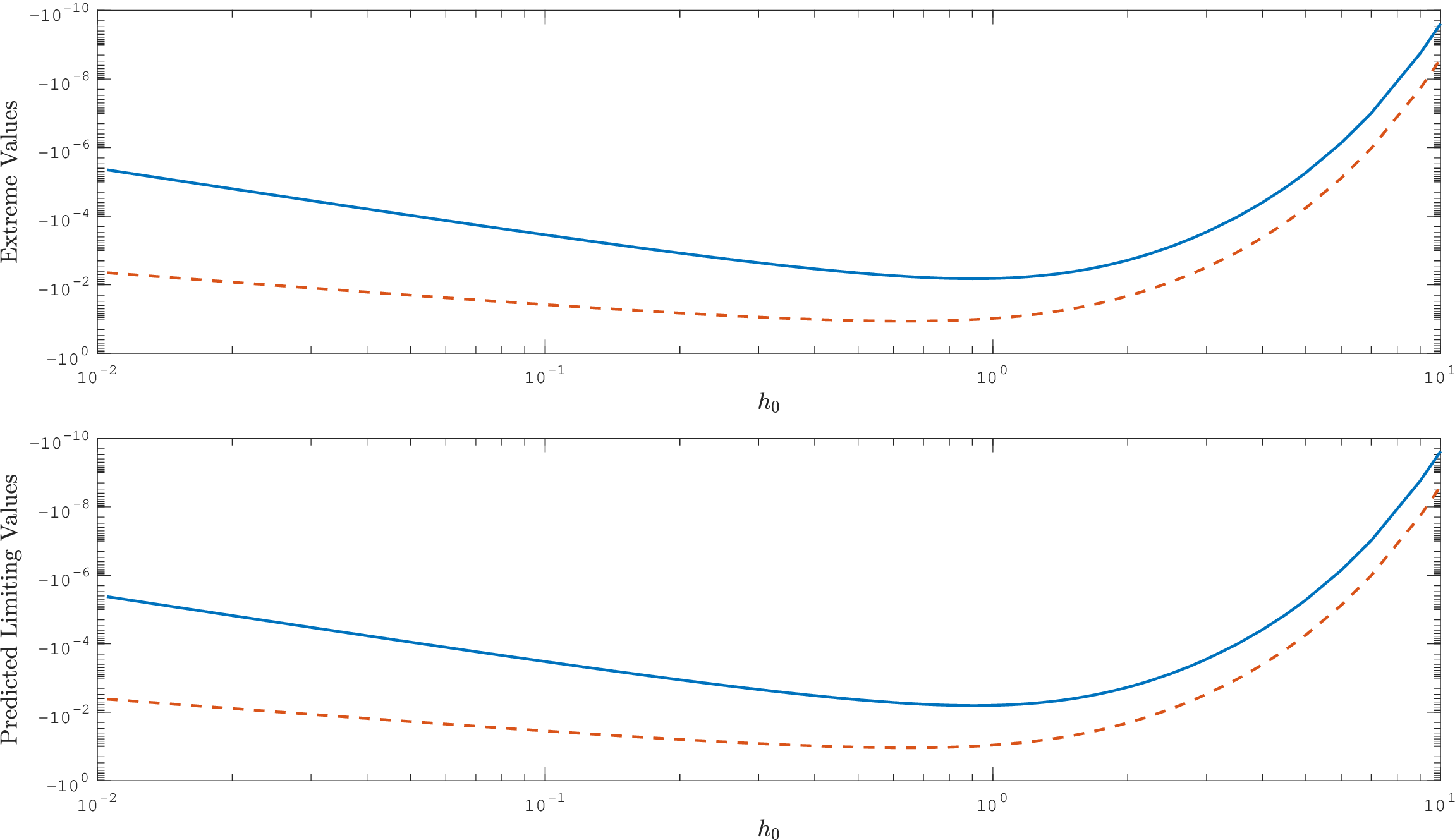}
    \caption{Change in the maximum (top) and limiting values (bottom) of the Bernoulli constant $B$ (solid blue line), inner product $\langle \eta,(\partial_hK_h)\eta\rangle/2\pi$ (dashed red line) as depth changes. Note that these are bounded away from 1.}
    \label{fig:extremeBernoulliEtaKhEta}
\end{figure}

Figure \ref{fig:extremeBernoulliEtaKhEta} shows the maximum and the limiting waves of the Bernoulli constant $B$ and  $\langle \eta,(\partial_h K_h)\eta\rangle$ divided by $2\pi$. We observe that these quantities do not fit into powers of the dispersion relation $\omega^2(kh)=kh\tanh(kh)$.

\section{Conclusion}
\label{sec-conclusion}

In this work, we have systematically developed the conformal formulation for Stokes waves on water of finite depth. 
We have incorporated two constraints, which follows from the Cauchy--Riemann equations, into the Lagrangian formulation 
via suitable Lagrange multipliers. This approach yields a closed system of constrained equations of motion in conformal variables 
and leads to Babenko's equation for traveling wave profiles supplemented by the corresponding integral constraint for Bernoulli's constant.
This conformal formulation is the main technical result of our work.

Within this framework, we carried out the existence and spectral stability analysis of the traveling wave solutions. We proved that co-periodic perturbations exchange stability precisely at the extremal points of the wave energy, or equivalently, at the extremal points of the horizontal momentum. The zero eigenvalue of the linearized operator has geometric multiplicity equal to $2$ and algebraic multiplicity equal to $4$ in the general case. However, 
the algebraic multiplicity is increased to at least $6$ at the extremal point of the wave energy. This result recovers the classical energy-momentum criterion in the finite-depth case \cite{Z-S86} and justifies it using a conformal mapping.

We have approximated the solution family numerically from the small-amplitude regime towards the steep Stokes waves. We discovered changes in the correspondence between the critical points of $\gamma$ and the kernel of the Jacobian $\mathcal{L}_h$. We point out 
that each critical point of the wave energy $\mathcal{H}$ occurs to the left 
of the critical point of $\gamma$ for a fixed conformal depth, in agreement with \cite{SB25}, and each critical point of the wave speed $c$ occurs to the right of the points where $\overline{\partial_h \eta} = 1$, which corresponds with the critical point of $\gamma$ for a fixed physical depth. The critical points of $\gamma$ for a fixed conformal depth $h$ discussed in \cite{SB25} precede those for a fixed physical depth $h_0$ discussed here. We conjecture that this sequence of critical points is repeated infinitely many times along the solution family towards the limiting Stokes wave of a peaked profile.

The formulation developed in this paper provides a solid basis for further analytical and numerical investigation of Stokes waves in finite depth. In particular, it opens the possibility of studying the modulational instability analytically and performing numerical continuation of solution branches, similar to \cite{BDS26} and \cite{DMP26}.


\begin{thebibliography}{99}
	
\bibitem{Babenko87}	K. I. Babenko, ``Some remarks on the theory of surface waves of finite amplitude", Sov. Math. Dokl. {\bf 35} (1987) 599--603.

 	\bibitem{BO82} T. B. Benjamin and P. J. Olver, ``Hamiltonian structure, symmetries and conservation laws for water waves", J. Fluid Mech. 
 {\bf 125} (1982) 137--185.
 
 \bibitem{BMV23} M. Berti, A. Maspero, and P. Ventura, ``Benjamin--Feir instability of Stokes Waves in finite depth", Arch. Rational Mech. Anal. {\bf 247}  (2023) 91.
 
 \bibitem{Buffoni1} B. Buffoni, ``Existence and conditional energetic stability of capillary-gravity solitary water waves by minimisation", 
 Arch. Ration. Mech. Anal. {\bf 173} (2004)  25--68.
 
  \bibitem{Buffoni2} B. Buffoni, ``Conditional energetic stability of gravity solitary waves in the presence of weak surface tension", 
 Topol. Methods Nonlinear Anal {\bf 25} (2005)  41--68.
 
\bibitem{BDT2000a} B. Buffoni, E. N. Dancer, and J. F. Toland, ``The regularity and local bifurcation of steady periodic water waves", Arch. Rational Mech. Anal. {\bf 152} (2000) 207--240.

\bibitem{BDT2000b} B. Buffoni, E. N. Dancer, and J. F. Toland, ``The sub-harmonic bifurcations of Stokes waves", Arch. Rational
 	Mech. Anal. {\bf 152} (2000) 241--271.
 	
\bibitem{BDS26} E. Byrnes, B. Deconinck, and A. Semenova, ``Stokes waves in water of finite depth", J. Nonlinear Waves {\bf 2} (2026) e7 (24 pages).
 
 \bibitem{CD23} R. P. Creedon and B. Deconinck, ``A high-order asymptotic analysis of the Benjamin--Feir instability spectrum in arbitrary depth", J. Fluid Mech. {\bf 956} (2023) A29.
 
 \bibitem{CDT22} R. P. Creedon, B. Deconinck, and O. Trichtchenko, ``Highfrequency instabilities of Stokes waves", J. Fluid Mech. {\bf 937} (2022) A24.
 
 \bibitem{DDS24}  B. Deconinck, S. Dyachenko, and A. Semenova, ``Self-similarity and recurrence in stability spectra of near	extreme Stokes waves", J. Fluid Mech. {\bf 995} (2024) A2.
 
 \bibitem{Mats} T. Dohnal, M. Ehrnstr\"{o}m, and J. U. Marstrander, ``A variational construction of solitary waves for the gravity Babenko equation 
 on finite depth", preprint (2026).

 \bibitem{DKSZ96} A. I. Dyachenko, E. A. Kuznetsov, M. D. Spector, and V. E. Zakharov, ``Analytical description of the free surface dynamics of an ideal fluid (canonical formalism and conformal mapping)", Phys. Lett. A
 {\bf 221} (1996) 73--79.

 \bibitem{DH19}  S. A. Dyachenko, and V. M. Hur. ''Stokes waves with constant vorticity: I. Numerical computation." Studies in Applied Mathematics {\bf 142.2} (2019): 162-189.


 \bibitem{DH26} S. A. Dyachenko, V. M. Hur, and D. A. Silantyev. ''Almost extreme waves." Journal of Fluid Mechanics {\bf 955} (2023): A17.
 
 \bibitem{DLK16} S. A. Dyachenko, P. M. Lushnikov, and A. O. Korotkevich. ``Branch cuts of Stokes wave on deep water.
 Part I: numerical solution and Pad\'e approximation", Stud. Appl. Math. {\bf 137} (2016) 419--472.
 
  \bibitem{DMP26} S. Dyachenko, R. Marangell, and D. E. Pelinovsky, ``Recurrent bifurcations of stability spectra for steep Stokes waves in a deep fluid", J. Nonlinear Waves {\bf 2} (2026) e17 (60 pages).
  
 \bibitem{DP25} S. Dyachenko and D. E. Pelinovsky, ``Bifurcations of unstable eigenvalues for Stokes waves derived from conserved energy", Physica D {\bf 483} (2025) 134925.
 
 \bibitem{DS23} S. A. Dyachenko and A. Semenova, ``Quasiperiodic perturbations 
 of Stokes waves: Secondary bifurcations and stability", 
 J. Comput. Physics {\bf 492} (2023) 112411 (13 pages)
 
 \bibitem{DS23SAPM} S. A. Dyachenko and A. Semenova, ``Canonical conformal variables based method for stability of Stokes waves", 
 Stud. Appl. Math. {\bf 150} (2023) 705--715.
 
 \bibitem{Ifrim1} B. Harrop--Griffiths, M. Ifrim, and D. Tataru, ``Finite depth gravity water waves in holomorphic coordinates", Ann. PDE {\bf 3} (2017) 1--102. 
 
 \bibitem{Ifrim2} J. K. Hunter, M. Ifrim, and D. Tataru, ``Two--dimensional water waves in holomorphic coordinates",  Commun. Math. Phys. {\bf 346} (2016) 483--552.

 \bibitem{KPS08} C. Kharif, E. Pelinovsky, and A. Slunyaev, {\em Rogue waves in the ocean} (Springer, New York, 2008)
 
\bibitem{LP25} S. Locke and D. E. Pelinovsky, ``On smooth and peaked traveling waves in a local model for shallow water waves", J. Fluid Mech. {\bf 1004} (2025) A1.
    
\bibitem{L-H 1} M. S. Longuet-Higgins, "The Instabilities of Gravity Waves of Finite Amplitude in Deep Water. I Superharmonics", Royal Society of London {\bf 360} (1978) 471-488.

\bibitem{L-H 2} M. S. Longuet-Higgins, "The Instabilities of Gravity Waves of Finite Amplitude in Deep Water. II Subharmonics", Royal Society of London {\bf 360} (1978) 489-505.

\bibitem{LHF1} M. S. Longuet-Higgins and M. J. H. Fox. "Theory of the almost-highest wave: the inner solution." Journal of Fluid Mechanics {\bf 80} (1977) 721--741.

\bibitem{LHF2} M. S. Longuet-Higgins and M. J. H. Fox. "Theory of the almost-highest wave. Part 2. Matching and analytic extension." Journal of Fluid Mechanics {\bf 85} (1978) 769--786.

\bibitem{NPW25} F. Natali, D.E. Pelinovsky, and S. Wang, ``Instability of the peaked traveling wave in a local model for shallow water waves", J. Nonlinear Waves {\bf 1} (2025) e10 (26 pages)

\bibitem{Saffman85} P.~G. Saffman, ``The superharmonic instability of finite-amplitude water waves'', J. Fluid Mech. {\bf 159} (1985) 169--174.

\bibitem{SB25} A. Semenova, and E. Byrnes, ``Stokes Waves in Finite Depth Fluids" Water Waves {\bf 7} (2025) 35--49.
	
\bibitem{Tanaka83} M.~Tanaka, ``The stability of steep gravity waves'', J. Phys. Soc. Japan {\bf 52} (1983) 3047--3055.
	
\bibitem{Tanaka85} M.~Tanaka, ``The stability of steep gravity waves. Part 2'', J. Fluid Mech. {\bf 156} (1985) 281--289.


\bibitem{Tanaka86} M.~Tanaka, "The stability of solitary waves", Phys. Fluids {\bf 29} (1986) 650–655.

\bibitem{Zakharov96} V. E. Zakharov and A. I. Dyachenko, ``High-Jacobian approximation in the free surface dynamics of an ideal fluid", Phys. D
{\bf 98} (1996) 652--664.
	
\bibitem{Z-S86} J.~A. Zufiria and P.~G. Saffman, "The Superharmonic Instability of Finite-Amplitude Surface Waves on Water of Finite Depth ", Elsevier {\bf74} (1986).
	


	

	

\end{thebibliography}
\end{document}